\documentclass{amsart}

\usepackage{etex}
\usepackage{amsmath, amssymb}
\usepackage{array}
\usepackage[frame,cmtip,arrow,matrix,line,graph,curve]{xy}
\usepackage{graphpap, color, paralist, pstricks}
\usepackage[mathscr]{eucal}
\usepackage[pdftex]{graphicx}
\usepackage[pdftex,colorlinks,backref=page,citecolor=blue]{hyperref}
\usepackage{pifont}
\usepackage{cancel}
\usepackage{mathtools}
\usepackage{tikz}

\usepackage{bbm}
\usepackage{dsfont}

\newcommand{\DD}{\mathbb{D}}

\newcommand{\pr}{\mathbb{P}}

\newcommand{\E}[0]{\mathbb{E}}

\newcommand{\beq}[1]{\begin{equation}\label{#1}}
\newcommand{\enq}[0]{\end{equation}}

\newcommand{\mn}[0]{\medskip\noindent}
\newcommand{\nin}[0]{\noindent}

\newcommand{\sub}[0]{\subseteq}
\newcommand{\sm}[0]{\setminus}
\renewcommand{\dots}[0]{,\ldots,}

\newcommand{\ov}[0]{\overline}

\newcommand{\less}[0]{~~\mbox{\raisebox{-.6ex}{$\stackrel{\textstyle{<}}{\sim}$}}~~}
\newcommand{\more}[0]{~~\mbox{\raisebox{-.9ex}{$\stackrel{\textstyle{>}}{\sim}$}}~~}
\newcommand{\lesser}[0]{~~\mbox{\raisebox{-.9ex}{$\stackrel{\textstyle{<}}{\approx}$}}~~}
\newcommand{\morer}[0]{~~\mbox{\raisebox{-.9ex}{$\stackrel{\textstyle{>}}{\approx}$}}~~}

\newcommand{\ra}[0]{\rightarrow}
\newcommand{\Ra}[0]{\Rightarrow}

\newcommand{\0}[0]{\emptyset}

\newcommand{\C}[2]{{{#1}\choose{{#2}}}}
\newcommand{\Cc}[0]{\tbinom}
\newcommand{\ga}[0]{\alpha }
\newcommand{\gb}[0]{\beta }
\newcommand{\gc}[0]{\gamma }
\newcommand{\gd}[0]{\delta }
\newcommand{\gD}[0]{\Delta }
\newcommand{\gG}[0]{\Gamma }

\newcommand{\gl}[0]{\lambda }

\newcommand{\go}[0]{\omega}
\newcommand{\gO}[0]{\Omega}

\newcommand{\gs}[0]{\sigma}
\newcommand{\gS}[0]{\Sigma}
\newcommand{\gz}[0]{\zeta}
\newcommand{\eps}[0]{\varepsilon }
\newcommand{\vt}[0]{\vartheta}
\newcommand{\vs}[0]{\varsigma}
\newcommand{\vr}[0]{\varrho}

\newcommand{\prh}[1][]{\pr_h}
\newcommand{\oH}[1][]{\ov{H}}
\newcommand{\oB}[1][]{\ov{B}}

\newtheorem{thm}{Theorem}[section]
\newtheorem{prop}[thm]{Proposition}

\newtheorem{lemma}[thm]{Lemma}
\newtheorem{obs}[thm]{Observation}
\newtheorem{cor}[thm]{Corollary}

\newtheorem{conj}[thm]{Conjecture}
\theoremstyle{definition}

\newtheorem{claim}[thm]{Claim}

\newcommand{\mbone}[0]{\mathbbm{1}}

\newcommand{\KKK}[0]{H}
\newcommand{\QQQ}[0]{A}
\newcommand{\kkk}[0]{k}
\newcommand{\uuu}[0]{u}

\begin{document}

\title{Asymptotics for palette sparsification}

\author{Jeff Kahn \and Charles Kenney}
\thanks{Department of Mathematics, Rutgers University}
\thanks{JK was supported by NSF Grant DMS1954035}
\email{jkahn@math.rutgers.edu,ctk47@math.rutgers.edu}
\address{Department of Mathematics, Rutgers University \\
Hill Center for the Mathematical Sciences \\
110 Frelinghuysen Rd.\\
Piscataway, NJ 08854-8019, USA}

\begin{abstract}

It is shown that the following holds for each $\eps >0$.
For $G$ an $n$-vertex graph of maximum degree $D$ 
and ``lists'' $L_v$ ($v\in V(G)$) chosen independently and uniformly from 
the ($(1+\eps)\ln n$)-subsets of $\{1\dots D+1\}$,
\[
\mbox{$G$ admits a proper coloring $\gs$ with $\gs_v\in L_v$ $\forall v$}
\]
with probability tending to 1 as $D\ra \infty$.

This is an asymptotically optimal version of a recent ``palette sparsification'' theorem of 
Assadi, Chen and Khanna.

\end{abstract}

\maketitle

\section{Introduction}\label{Intro}

Recall that a \emph{list coloring} of a graph $G$ (on vertex set $V$)
is a proper vertex coloring in which each vertex receives a color from some
prescribed list (particular to that vertex).  Formally, 
for $L=(L_v:v\in V)$ with $L_v\sub \gG = \{\mbox{``colors''}\}$,
an \emph{L-coloring} of $G$ is a (proper) coloring $\gs$ of $V$ with
$\gs_v\in L_v$ $\forall v$.  If such exists then we also say $G$ is \emph{L-colorable}.
(See \cite{Vizing,ERT} for origins and e.g.\ \cite{Diestel} for general list coloring background.)

\mn

This work was motivated by the following ``palette sparsification'' 
theorem of Assadi, Chen and Khanna~\cite{ACK};
see \cite{ACK} for the surprising algorithmic consequences.
(Here and throughout the paper we use $n=|V|$ and $\log=\ln$.)
\begin{thm}\label{TACK}
There is a fixed C such that if $\gD_G = D$ and the $L_v$'s are drawn uniformly and independently
from the ($C\log n$)-subsets of
$\gG=[D+1]$ $(:=\{1\dots D+1\})$, then G is L-colorable w.h.p.\footnote{\emph{with high probability},
meaning with probability tending to 1 as $n\ra \infty$.} 
\end{thm}

\nin
(Strictly speaking, the $C\log n$ should be $\min\{C \log n, D+1\}$, with the statement trivial if 
this is $D+1$.)

Notice that Theorem~\ref{TACK} fails when $C=1$.  For example, if $G=K_{D+1}$, then an $L$-coloring
is the same thing as a perfect matching in the bigraph on $(V,\gG)$ with $v\sim \gc$ iff $\gc\in L_v$;
and when $C=1$ a simple (and standard) calculation says that
$\pr(\cup_v L_v = \gG)$ is only about $e^{-1}$.
Here we prove the following natural guess, which we first heard suggested (enthusiastically)
in a lecture of Sepehr Assadi.   
\begin{thm}\label{TKK}
For any fixed $\eps >0$ the following is true.
If $\gD_G \leq D\leq n$ and the $L_v$'s are drawn uniformly and independently
from the ($(1+\eps)\log n$)-subsets of
$\gG=[D+1]$, then G is L-colorable w.h.p.
\end{thm}

\nin
Equivalently, for $\vs$ a slow enough $o(1)$,
the theorem holds with $\eps$ replaced by $\vs$.
(In fact any $\vs =\go(1/\log n)$, which is what's needed when $G $ is $K_{D+1}$, should 
suffice, but we don't show this.)

\mn

While coloring with random sublists of given lists is (\emph{in retrospect}) natural
in the algorithmic setting of~\cite{ACK}, it is also of considerable interest on purely combinatorial grounds.
We will say a little more about this, and about possible extensions of Theorem~\ref{TKK}, in
Section~\ref{Disc}.

The proof of Theorem~\ref{TKK}, which turned out to be (much) trickier
than expected, is previewed in Section~\ref{Sketch}, at the end of which we will be in 
a better position to say something intelligible about how the rest of paper is organized.

\section{Sketch}\label{Sketch}

Here we give a quick overview of the proof of Theorem~\ref{TKK}, 
noting that while
the strategies described are reasonably natural, 
their implementation is delicate and not very straightforward
(and that most of what's sketched here will be easier to say comprehensibly in the more formal settings 
to come).

\mn

It will be convenient to assume in what follows $G$ is $D$-regular, which is justified by the next observation.
\begin{prop}\label{Preg}
Any $n$-vertex (simple) graph $G$ with maximum degree at most $D$ is contained in a $D$-regular
(simple) graph $G'$ on $n'\leq n+D+2$ vertices.
\end{prop}
\nin
Here we just need $n'=n^{1+o(1)}$ to make $\log n'\sim \log n$.
We expect Proposition~\ref{Preg} is known but couldn't find a reference with such a 
bound on $n'$, so will sketch a proof at the end of this section.

\mn

We repeat the assignment:
We are given a fixed $\gd>0$ and should show that for 
a $D$-regular $G$ and $L=(L_v:v\in V)$
as in Theorem~\ref{TKK}---namely, the $L_v$'s are independent random lists of size 
\beq{ell}
\kkk:=(1+\gd)\log n 
\enq
---$G$ is $L$-colorable w.h.p.

As in \cite{ACK} we begin by partitioning $G$ into a ``sparse'' part, $V^*$, and
small, quite dense ``clusters'' (called ``almost-cliques'' in \cite{ACK}).
We will use the following statement, extracted from \cite[Lemma 2.3]{ACK} 
(itself a variant of a decomposition from \cite{HSS}).
\begin{lemma}\label{LACK}
For any D, $\eps>0$, and $D$-regular graph $G$ on vertex set $V$, 
there is a partition
$
V=V^*\cup C_1\cdots\cup C_m
$
such that:

\nin
{\rm (a)}  
$~v\in V^* ~\Ra ~ |\{w\sim v: d(v,w) < (1-\eps)D\}| > \eps D$;

\nin
{\rm (b)} 
$\,$ for all $i\in [m]$:  

{\rm (i)} \, $|C_i|\in [(1-\eps)D, (1+6\eps)D]$, and

{\rm (ii)}\, for all $v\in C_i$, $|N_v\sm C_i|< 7\eps D$ and $|C_i\sm N_v|< 6\eps D$.
\end{lemma}
\nin
(As one would expect, the particular constants are unimportant;  we will drop them 
later, but keep them here to show the connection with \cite{ACK}.)
An easy consequence of (a) (which is all we will use it for) is
\beq{generic'} 
|\ov{G}[N_v]| > (\eps^2/2)D^2\,\,\,\,\forall v\in V^*.
\enq

We will use Lemma~\ref{LACK} with $\eps$ a small constant (depending on $\gd$),
to which we will not assign a specific value; constraints on what that value might be
will appear mostly in the ``dense'' phase of the argument, beginning in Section~\ref{Clusters}
(see ``Note on parameters'' near the beginning of that section).

\mn

As in \cite{ACK} we color in two stages, the first (``sparse'') coloring $V^*$ and the second 
(``dense'') extending to the $C_i$'s.

\nin
\textbf{Sparse phase.}  Our aim here
is to reduce to the following result of 
Reed and Sudakov \cite{RS}.
For  $v\in V$ and $\gc\in \gG$, define the ``color degree'' 
\[
d_\gc(v):=|\{w\sim v:\gc\in L_w\}|.
\]
\begin{thm}\label{TRS}
For any $\vr>0$ there is a $t_\vr$ such that
if $t> t_\vr$, 
all color degrees are at most $t$, and all list sizes are at least $(1+\vr)t$, 
then $G$ is $L$-colorable.
\end{thm}
\nin
\emph{Remarks.}
(a)  Equivalently, this is true with $\vr=\vr(t)$ a sufficiently slow $o(1)$.
Reed \cite{Reed} conjectured that  
Theorem~\ref{TRS} is true if the list sizes are at least $t+1$; this was disproved---but 
not with $t+1$ replaced by $t+2$---by Bohman and Holzman \cite{BH}.  
As observed in both these references, it would be interesting to see if the extra
$\vr t$ in Theorem~\ref{TRS} could be reduced, perhaps even to $O(1)$.

\nin
(b)  Theorem~\ref{TRS}, a natural tool in the present context, is also used in \cite{AA},
though in a situation where its relevance is more immediate than will be the case here.

\nin

We will use Theorem~\ref{TRS} with lists that are subsets of the $L_v$'s of Theorem~\ref{TKK}.
Note that for the $L_v$'s themselves the typical color degree, 
$D \kkk/(D+1)$, is essentially equal to list size.  

Our initial plan here (which doesn't entirely work) is as follows.
We first tentatively assign single random colors to the vertices 
(the color assigned to $v$ being a first installment on $L_v$), retaining these 
colors at the properly colored vertices
(those assigned colors not assigned to their neighbors) and discarding the rest. 
We then choose the balances of the lists and, for $v$ not colored in the first stage,
confine ourselves to the surviving colors in $L_v$ (those not already used at neighbors of $v$).
This is supposed to work as follows.

For a given $v$ the first stage above should yield:

\nin
(a)  a positive fraction (essentially $e^{-1}$) of $v$'s neighbors properly colored;

\nin
(b)  (\emph{because we are in the sparse phase)}
significant duplication among the colors assigned
to these neighbors, so that the number of such colors 
is significantly less 
than the number of colored neighbors.

\nin
On the other hand, it is easy to see that if (a) and (b) (suitably quantified) hold for all $v$, then 
the \emph{natural} values (with respect to still-allowed colors)
of the list sizes and color degrees support use of Theorem~\ref{TRS}.

But of course we need to say that w.h.p.\ all 
(or a strong enough \emph{most}) of these
desirable behaviors (i.e.\ (a), (b) and the ``natural'' values above)
happen \emph{simultaneously}.
Modulo some adjustments this is (provably) true for $n$ less than some $e^{\Theta(D)}$,
but not all the way to $e^D$ (beyond which $\ell >D$ and there is nothing to show).
This was a surprise, since larger $n$, with its larger lists, ought to be easier, and
the way we wound up handling it seems to us one of the more interesting features of the present 
analysis; but we postpone further discussion of this to Section~\ref{Generic}.

\nin
\textbf{Dense phase.}
We process the $C_i$'s one at a time, in each case assuming an arbitrary coloring
(even an arbitrary \emph{labeling} would be okay here) of $V^*$ and the earlier clusters.

For a given $C_i$ we may take as our point of departure the bipartite matching discussion
following Theorem~\ref{TACK}---where we want to match vertices of $C_i$ with colors 
from their lists, so expect some use of Hall's Theorem.  There are two complications here:

First, we can have $|C_i|>|\gG|$, so our 
``matching'' must actually use some colors more than once.  
Second, whatever the size of $C_i$, our situation differs from the one following
Theorem~\ref{TACK} in that a vertex $v$ can only be matched with colors from
$L_v$ that haven't already been used in the neighborhood of $v$.

A first ingredient in coping with these issues (in more interesting cases) 
is use of a small portion of $\gG$ to 
assign common colors to enough pairs of nonadjacent
vertices to reduce to a situation where the number of available colors does somewhat
exceed the number of uncolored vertices.
(See ``Process'' in 
Section~\ref{Clusters}.
Something of this type, though under less severe constraints, is also used in \cite{ACK},
and seems more or less unavoidable.)
Even after this, the ``Hall'' part of the argument---showing good matching
behavior in a random bigraph with some disallowed edges---is less straightforward
than might be expected, the most interesting point here being Lemma~\ref{Lswitch}
(and the need for it).

After reviewing a few standard tools in Section~\ref{Prelim},
we give the basic argument for the sparse phase in Section~\ref{Generic},
postponing the proofs of Lemmas~\ref{LPh1a} and \ref{LC0}, 
which carry most of the load.
The martingale concentration machinery needed for these is reviewed in Section~\ref{Mart},
and the lemmas themselves are proved in Section~\ref{Conc}.
The dense phase is treated in Sections~\ref{Clusters}-\ref{Matchings}, again
with the main line of argument in Section~\ref{Clusters}, followed by proofs
of the supporting Lemmas~\ref{Lsucceed}-\ref{Lswitch}
in Sections~\ref{PLsucceed} and \ref{Matchings}.
(The sparse and dense parts can be read independently.)

\nin
\textbf{Usage}

As above $V=V(G)$ is of size $n$ and $\gG = [D+1]$.
Any $\gc$, perhaps subscripted, is by default a member of $\gG$.
In graphs, $\nabla(X,Y)$ is the set of edges joining the disjoint vertex sets $X,Y$;
and $\nabla(X)=\nabla(X,V\sm X)$, $N(\cdot)$ is neighborhood,
and $N(S)=\cup_{v\in S}N(v)$; and $H[X]$ is the subgraph of $H$ induced by $X$.

Asymptotics:  we use $o(\cdot)$, $\go(\cdot)$, $\sim$, etc.\ in the usual ways, and often write 
$a\less b$ for $a<(1+o(1))b$  (all as $D\ra \infty$).
The discussion, especially in the dense phase, will involve various small and large constants,
but the implied constants in $O(\cdot)$, $\gO(\cdot)$, $\Theta(\cdot)$
do not depend on these.

We will always assume $D$ is large enough to support
our assertions and, following a common abuse, usually---namely, unless it makes a difference---pretend 
large numbers are integers.
In most of what follows, the size of $\gG$ will appear far more often than the degree,
so to avoid cluttering the paper with $(D+1)$'s we proceed as follows.
For the sparse phase, the difference is irrelevant (the argument doesn't really
change if we shrink $|\gG|$ by 1), and we ignore it, using $D$ for both degree and $|\gG|$.
But when we come to the dense phase, the difference does sometimes matter (cf.\ the paragraph
following Theorem~\ref{TACK}), so in the last few sections (beginning with Section~\ref{Clusters})
we occasionally (sparingly) 
use $\DD$ for $D+1=|\gG|$.

We will sometimes use square brackets to indicate quick comments and justifications.

\mn

\nin
\emph{Proof of Proposition~\ref{Preg}} (sketch).  
We may assume (by adding edges in $V=[n]$ until each edge of $\ov{G}$ has an 
end with $G$-degree $D$)
that $W:=\{v\in V: d_G(v)<D\}$ supports a clique of $G$ (so $|W|\leq D$).
We claim that we may obtain a $G'$ as in the proposition by doing:

\begin{itemize}

\item
let $X$ be a set of vertices disjoint from $V$, with 
$|X|\in \{D+1,D+2\}$ of the same parity as $n$;

\item
choose $K$ bipartite on $W\cup X$ with $d_K(v) = D-d_G(v)$ $\forall v\in W$ and
$|d_K(x)-d_K(x')|\leq 1$ $\forall x,x'\in X$; 

\item
let $F$ be a graph on $X$ with $d_F(x) = D-d_K(x)$ $\forall x\in X$;

\item
and let $G' = G\cup K\cup F$.

\end{itemize}

\nin
If achievable, $G'$ clearly has the desired properties, and the only questionable item
is existence of $F$.  This will be given by the next observation, which 
is easily proved by induction using the Havel-Hakimi Theorem (\cite{Havel,Hakimi} or e.g.\
\cite[Thm.\ 1.3.31]{West}).
\begin{prop}
A nonnegative integer sequence $(f_1\dots f_m)$ with $|f_i-f_j|\leq 1$ $\forall i,j$
is the degree sequence of a simple graph iff
(i) $m\geq \max f_i +1$ and (ii) $\sum f_i$ is even.
\end{prop}
\nin
Here we just need (ii), which is again easy:  setting $s=\sum_{v\in W}(D-d_G(v))$ and 
$f_x = D-d_K(x)$ (for $x\in X$), we have (with $\equiv$ denoting congruence mod 2)
\[
\sum_{x\in X}f_x = |X|D-s \equiv nD-s = 2|E(G)|.
\]\qed

\section{Preliminaries}\label{Prelim}

We first recall, presumably unnecessarily, the Lov\'asz Local Lemma (\cite{ErLov} and e.g.\
\cite[Ch.\ 5]{AS} for more general version):
\begin{lemma}
Suppose $A_1,\ldots, A_m$ are events in a probability space with
$\pr(A_i)\leq p$ $\forall i$,
and $\gG$ is a graph on $[m]$
with $A_i$ independent of $\{A_j: j\not\sim_\gG i\}$ ($\forall i$). 
If $ep(\gD_\gG +1)<1$, then
$
\pr(\cap \bar{A_i}) >0.
$
\end{lemma}

Most of what we do below relies on various
concentration assertions, beginning with the following
standard
``Chernoff"
bounds.
Recall that a random variable $\xi$ is \emph{hypergeometric} if, for some
$s,a$ and $k$,
it is distributed as $|X\cap A|$, where $A$ is a fixed $a$-subset
of the $s$-set $S$ and $X$ is uniform from $\C{S}{k}$.
\begin{thm}
\label{T2.1}
If $\xi $ is binomial or hypergeometric with  $\mathbb{E} \xi  = \mu $, then for $t \geq 0$,
\begin{align}
\pr(\xi  \geq \mu + t) &\leq
\exp\left[-\mu\varphi(t/\mu)\right] \leq
\exp\left[-t^2/(2(\mu+t/3))\right], \label{eq:ChernoffUpper}\\
\pr(\xi  \leq \mu - t) &\leq
\exp[-\mu\varphi(-t/\mu)] \leq
\exp[-t^2/(2\mu)],\label{eq:ChernoffLower}
\end{align}
where $\varphi(x) = (1+x)\log(1+x)-x$
for $ x > -1$ and $ \varphi(-1)=1$.
\end{thm}
\nin
(See e.g.\ \cite[Theorems 2.1 and 2.10]{JLR}.)
For larger deviations the following consequence of the finer bound in \eqref{eq:ChernoffUpper}
is helpful.  
\begin{thm}
\label{Cher'}
For $\xi $ and $\mu$ as in Theorem~\ref{T2.1} and any $K$,
\begin{eqnarray*}
\pr(\xi  > K\mu) < \exp[-K\mu \log (K/e)].
\end{eqnarray*}
\end{thm}
\nin

We defer to Section~\ref{Mart} the more sophisticated martingale concentration
statements that will support our treatment of 
the ``sparse'' phase in Section~\ref{Conc}.

We also need the simplest case of the Janson Inequalties
(see \cite{Janson} or \cite{JLR}, in particular
\cite[Theorem 2.18(ii)]{JLR} for the case in question).  
Suppose $A_1\dots A_m$ are subsets of a
finite set $S$; let $E_i$ be the event $\{S_p\supseteq A_i\}$
(with $S_p$ the usual ``$p$-random'' subset of $S$);
and set $\mu =\sum\pr(E_j)$ and
\[   
\ov{\gD} = \sum\sum\{\pr(E_iE_j): A_i\cap A_j\neq\0\}.
\]   
(Note this includes
diagonal terms.)

\begin{thm}\label{TJanson}
With notation as above, 
$\,\,\pr(\cap \ov{E}_i) \leq \exp[-\mu^2/\ov{\gD}]$.
\end{thm}

\section{Coloring $V^*$}\label{Generic}

Recall (again) that we are trying to color $G$ from random lists $L_v$ of size $\kkk=(1+\gd)\log n$
(see \eqref{ell}).

The goal of the present section is to produce an $L$-coloring of $V^*$.
For this discussion we retain the cluster vertices since it will be convenient to work
with a regular graph;
but list choices for cluster vertices in this phase are ``dummies'' to be discarded once they
have served their purpose here.

Fix $\eps$ small enough (depending on $\gd$) to support what follows.
(As noted following Lemma~\ref{LACK}, what this entails will appear mainly in 
Section~\ref{Clusters}; again, see ``Note on parameters'' following \eqref{notext}.)

  For the rest of the paper
assume the decomposition $V=V^*\cup C_1\cdots\cup C_k$
of Lemma~\ref{LACK}, and, 
setting $\vt =\eps^2/2$, recall \eqref{generic'}, which is all we will need from 
the decomposition in the sparse phase:
\beq{generic} 
|\ov{G}[N_v]| > \vt D^2\,\,\,\,\forall v\in V^*.
\enq

As suggested in Section~\ref{Sketch} (see (a), (b) under ``Sparse Phase''),
we first want to say 
that we can (w.h.p.) specify, at some affordably small cost in list size, 
$T\sub V$ and an $L$-coloring $\gs$ 
of $T$ such that
for every $v\in V^*$ (though we only care about $v\not\in T$),
\beq{ideal1}
|T\cap N_v|\sim e^{-1}D
\enq
and, with $\vt'$ to be specified below (see \eqref{c'}),
\beq{ideal2}
|T\cap N_v| - |\gs(T\cap N_v)| > \vt' D.
\enq
(Membership in $V^*$ is irrelevant for \eqref{ideal1}, which we
will really show for \emph{all} $v$, 
but the restriction  is, obviously,
needed for \eqref{ideal2}.)

A natural way to try to arrange this is as follows.  Choose
\beq{tau1}
\mbox{$\tau $ uniformly from $\gG^V$}
\enq
(taking $\tau_v$ to be a first member of $L_v$)
and let
\beq{Tgs}
\mbox{$T=\{v:\tau_w\neq\tau_v~\forall w\sim v\}\,\,$ and
$\,\,\gs   =\tau|_T$.}
\enq
Thus we assign tentative colors $\tau_v$, uncolor wherever $\tau$ is improper,
and let $T$, $\gs$ be what's left.

Before proceeding we check that \eqref{ideal1} and \eqref{ideal2} are correct at the 
level of expectations.  For \eqref{ideal1} this is true---that is,
\beq{ideal1'}
\E |T\cap N_v|\sim e^{-1}D
\enq
---since 
\beq{przT'}
\pr(z\in T) = (1-1/(D+1))^D\sim e^{-1}.
\enq

For \eqref{ideal2} the story is a little longer (and we do need $v\in V^*$,
or, more precisely, \eqref{generic}).  
It will be easier to work with the following lower bound.
For distinct $w,z\in N_v$ define the event
\[
A_{wz}^v =\{w,z\in T, \tau_w = \tau_z\neq \tau_x\,\,\forall x\in N_v\sm \{w,z\}\}
\]
(that is, $\{w,z\}\sub T\cap N_v$ is monochromatic [under $\tau$] and 
not contained in any monochromatic triple in $T\cap N_v$;
note this implies $w\not\sim z$).
Then
\[   
|T\cap N_v| - |\gs(T\cap N_v)|\geq
|\{\{w,z\}\sub N_v:\mbox{$A^v_{wz}$ holds}\}|.
\]   
We will show \eqref{ideal2} with its l.h.s.\ replaced by this lower bound---that is,
\beq{ideal2'}
|\{\{w,z\}\sub N_v:\mbox{$A^v_{wz}$ holds}\}| > \vt' D
\enq
---so now
think of the expectation version of \eqref{ideal2'}.

For $w,z$ nonadjacent neighbors of $v$, let (recalling that $N(S)=\cup_{z\in S}N(z)$)
\beq{Jvwz}
J^v_{wz} = N(v,w,z)\sm \{w,z\}.
\enq
We may rewrite 
\[
A^v_{wz}=\{\tau_w=\tau_z\neq \tau_x \,\forall x\in J^v_{wz}\},
\]
and observe that
\beq{prAvwz}
\pr(A^v_{wz})=D^{-1}(1-1/D)^{|J^v_{wz}|} > D^{-1}e^{-3}
\enq
(where, as promised at the end of Section~\ref{Sketch}, we now start pretending
$|\gG|=D$),
whence, in view of 
\eqref{generic}, 
\beq{Ewz}
\E |\{\{w,z\}\sub N_v:\mbox{$A^v_{wz}$ holds}\}|
> \vt De^{-3};
\enq
so in \eqref{ideal2} and \eqref{ideal2'} we take 
\beq{c'}
\vt' = 0.5 \vt  e^{-3} \,\,\,(= 0.25\eps^2 e^{-3}).
\enq

\mn

Given the above expectation versions 
we would
like to say, using martingale concentration and a union bound over $v$, that 
w.h.p.\ all left-hand sides in \eqref{ideal1} and \eqref{ideal2'} 
are close to their expectations;  which---again, as in Section~\ref{Sketch}---would give
some hope that the \emph{as-yet-unexamined} lists $L_v\sm \{\tau_v\}$ support 
application of Theorem~\ref{TRS}.
This turns out to work for $n$ up to some $\exp[\Theta(D)]$ (whereas 
we may assume $(1+\gd)\log n< D+1$), but does eventually break down because the concentration is not 
strong enough.  Handling this seemingly minor issue turned out to be (to date)
surprisingly challenging, and, as mentioned earlier,
the fix below seems to us one of the more interesting aspects of the present
work.

We will use said fix only for quite large $n$, say  
\beq{nlarge}
n > \exp[D^{.9}].
\enq
(It can be made to work in general, but this requires some care and, not needing it,
we won't bother.)

Call $n$ \emph{large} (\emph{small}) according to whether it satisfies (doesn't satisfy)
\eqref{nlarge}.  In Section~\ref{P23} we will show
\beq{small23}
\mbox{\emph{if $n$ is small then w.h.p.\ the T and $\gs$ of \eqref{tau1}-\eqref{Tgs}
satisfy \eqref{ideal1} and \eqref{ideal2'} for all $v\in V^*$.}}
\enq

For larger $n$ we proceed as follows.
Set $\ell =0.1\gd \log n $.  We choose the $\tau_v$'s in two steps:

\mn
\emph{Step 1}:
Each $v$ chooses (independently, uniformly) $L_v^0\in \C{\gG}{\ell}$, and we set
$L^0 = (L_v^0:v\in V)$.

\nin
\emph{Step 2}:  Each $v$ chooses $\tau_v $ (independently, uniformly) from $L_v^0$,
and we define $T,\gs$ as in \eqref{tau1} and \eqref{Tgs}.

\mn

The crucial gain here is that we now need w.h.p.\ statements \emph{only for the} $L_v^0$'s:
it is enough to say that w.h.p.\ these partial lists support
\emph{existence}
of a good 
$T$ and $\gs$:  
\beq{CL0}
\mbox{\emph{if n is large, then w.h.p.\ 
$L^0$ is
such that w.p.p.\ the $T,\gs$ of 
Step 2 satisfy \eqref{ideal1} and \eqref{ideal2'} for all $v \in V^*$},}
\enq   

\nin
where \emph{w.p.p.---with positive probability}---just means the probability is not zero.
Thus the 2-step procedure succeeds because concentration for the relevant quantities output by
Step 1---those in Lemma~\ref{LPh1a}---\emph{is} strong enough to support union bounds, and the somewhat 
weaker concentrations of Step 2 are (far) more than enough for use with the 
Local Lemma.

For coloring
$V^*\sm T$, we then forget what we know about $L^0$
and work with fresh,
slightly smaller lists---thus the cost in list size mentioned following \eqref{generic}---so as to 
avoid the no longer comprehensible law of $L^0$.

\mn

In the rest of this section we set up the proof of \eqref{CL0}---namely, specifying, in 
Lemma~\ref{LPh1a}, the (w.h.p.) properties of $L^0$ that will, as asserted in Lemma~\ref{LC0},
support the ``w.p.p.'' statement of \eqref{CL0}---and 
then finish the sparse phase assuming these lemmas, whose proofs are postponed to
Section~\ref{Conc}.
So for the next few paragraphs (through the statement of Lemma~\ref{LC0}) we 
continue to think of large $n$.

We now use $\pr'$ and $\E' $ for probabilities and expectations associated with 
Step 1, and $\pr''$ for probabilities associated with Step 2;
thus the $\pr''(\cdots)$'s are r.v.s determined by $L^0$.
We continue to use $\pr$ and $\E$ for probabilities and expectations
associated with the original experiment
(that is, $\tau_v$'s chosen uniformly and independently from $\gG$), observing that
\beq{EEP}
\mbox{$\E'\pr''(\cdot) =\pr(\cdot)~$ and $~\E'\E''[\cdot] =\E[\cdot].$}
\enq

\nin

Let $\theta\ra 0$ and $C\ra \infty$, both slowly , with
\beq{Ceps}
C^{\theta} =\go(1).
\enq
[The other (mild) constraints on $(C,\theta)$ are:  
$C\theta^{-2}=o(\ell)$ (see following \eqref{X2X}),
$C^{-\theta D}=o(\theta)$ in \eqref{EU}, and
$C = D^{o(1)}$ at the ends of the two proofs of \eqref{ideal1}
in Section~\ref{P23}; e.g.\ $\eps =\go(1/\log C)$, with $C\ra\infty$ slowly enough will do.]

\begin{lemma}\label{LPh1a}
W.h.p.\ (w.r.t.\ the choice of $L^0$)
\beq{whp1} 
\mbox{for all $z$,   $\,\,\,\pr''(z\in T)\sim \pr(z\in T) $;}
\enq
\beq{whp2}
\mbox{for all $v$,   
$\,\,\,|\sum_{\{w,z\}\sub N_v}[\pr''(A^v_{wz})- \pr(A^v_{wz}) ]|=o(D)$;}
\enq
\beq{whp3}
\mbox{for all $v$, $Y\in \C{N_v}{\theta D}$ and $\gG_0\in \C{\gG}{\theta D/C}$,
$\,\,\,\sum_{y\in Y} |L_y^0\cap \gG_0| < 2 \theta ^2 D\ell/C~$;}
\enq
and
for all $v$ and $\tau:N_v\ra \gG$,
\beq{whp4}
\sum_{w\not\sim v}|\tau(N_w\cap N_v)\cap L_w^0| < 2D\ell.
\enq
\end{lemma}

\begin{lemma}\label{LC0}
If the lists $L_v^0$ are as in Lemma~\ref{LPh1a}, then w.p.p.\
Step 2 gives $T,\gs$ as in \eqref{ideal1}, \eqref{ideal2}.
\end{lemma}

Again, the proofs of Lemmas~\ref{LPh1a} and \ref{LC0}
are deferred to Section~\ref{Conc}, and in the rest of this section we 
dispose of the sparse phase assuming the lemmas are correct.
Thus we assume we have chosen $\tau$ (preceded by $L^0$ if $n$ is large) with
associated $T$ and $\gs$ satisfying \eqref{ideal1} and \eqref{ideal2},
and show
that w.h.p.\ $\gs$ can be extended to $V':=V^*\sm T$.

For large $n$, ``w.h.p.'' here refers to the choice of ($L_v\sm L_v^0$)'s
(for $v\in V'$, lists for other $v$'s being no longer relevant); but it will be easier to 
choose lists $L_v'$ uniformly (and 
independently) from $\C{\gG}{t}$, where
\beq{tgd'}
t:=\kkk-\ell = (1+0.9\gd)\log n=:(1+\gd')\log n
\enq
(recall---see \eqref{ell}---$\kkk =|L_v| = (1+\gd)\log n$),
and show that \emph{these} give the desired extension w.h.p.
That the $L_v$'s also do this then follows, 
since
we can couple so $L_v'\sub L_v$.
Similarly, for small $n$---though this is very silly---we can take $t =\kkk-1$,
and again work with $L_v'$ uniform from $\C{\gG}{t}$.

By \eqref{ideal1} and \eqref{ideal2} 
the lists of still-allowed colors,
\[
L_v'':=L_v'\sm \gs(T\cap N_v),
\]
and color degrees, $d''_\gc(v) :=|\{w\in N'_v:\gc\in L_w''\}|$ (where $N'$ is
neighborhood in $V'$), satisfy
\[
\E |L_v''| \more (1-e^{-1} + \vt')t
\]
(recall $\vt'$ was defined in \eqref{c'}) and
\[
\E d''_\gc(v) \less (1-e^{-1} )t.
\]

Of course if \emph{actual} behavior agreed with these expectations then Theorem~\ref{TRS} 
would give the desired extension of $\gs$ to $V'$;
but in reality there will usually be some
exceptions, and we can invoke the theorem only after disposing of these, as follows.

We first observe that w.h.p.
\beq{bad1}
|L_v''|>\eta \log n \,\,\,\,\,\,\forall v \in V'
\enq
for a suitable fixed (positive) $\eta$.

\begin{proof}
Since $|\gs(T\cap N_v)|< D/e$ (as follows from \eqref{ideal1} and \eqref{ideal2}), 
\[
\pr(|L_v''|< \eta\log n) < \C{(1+\gd')\log n}{\eta\log n}e^{- (1+\gd'-\eta)\log n},
\]
which is $o(n^{-1})$ if $\eta\log (1/\eta)$ is somewhat less
than $ \gd'$.
\end{proof}

In choosing lists $L_v'$ (again, for $v\in V'$), we first choose the cardinalities
$
m_v:= |L_v'\cap \gs(T\cap N_v)|;
$
say $v$ is \emph{bad} if
\[
m_v> (e^{-1} - \vt'+\vs)t,
\]
with $\vs=\vt'/3$; and 
let $B$ be the set of bad vertices.
Then Theorem~\ref{T2.1} gives
\[
\pr(v\in B) < n^{-\vs' },
\]
with $\vs'=\Theta( \vs^2)$, implying that w.h.p.\ (\emph{very} conservatively when $n$ is 
large)
\beq{bad2}
|B\cap N_v| < 2D^{1-\vs' } \,\,\,\,\,\,\forall v.
\enq
[With $X=|B\cap N_v|$ and $Y\sim \textrm{Bin}(D,n^{-\vs'})$, 
\[
\pr(X\geq 2D^{1-\vs'})\leq\pr(Y\geq 2D^{1-\vs'}) =o(n^{-1})
\]
is a routine consequence of one of Theorems~\ref{T2.1}, \ref{Cher'} (say the latter
unless $n=O(D)$).]

\mn

Given the $m_v$'s, we choose the lists
$L_v''$, first for $v\in B$ 
and then for the remaining $v$'s.  
(So each $L_v''$ is uniform from the $(t -m_v)$-subsets of 
$\gG\sm \gs(T\cap N_v)$.
There is of course no point in specifying any $L_v'\sm L_v''$.)

We next choose an $L''$-coloring $\gs_B$ of $B$, noting that
a greedy coloring 
(it's best to think of choosing $L_v''$ when we are ready to color $v$)
produces such a $\gs_B$ w.h.p.
[Assuming the w.h.p.\ statement \eqref{bad1}---which depends only on $m_v$---and
using \eqref{bad2}, the probability that the process gets stuck at $v$ is
less than
\[
\left(|B\cap N_v|/D\right)^{|L_v''|} < (2D^{-\vs'})^{\eta\log n} =o( n^{-1}).]
\]

We then want to use Theorem~\ref{TRS} to color $V'':= V'\sm B$.
Thinking of choosing $L_v''$ for $v\in V''$ after we have specified $\gs_B$,
we next observe that
\eqref{bad2} implies that w.h.p.\ 
\beq{Lvtau}
|L_v''\cap \gs_B(B\cap N_v)| < (2/\vs')\log n/\log D =:J = o(\log n)\,\,\,\,\,\,\forall v\in V''.
\enq
[Given \eqref{bad2}, the probability that \eqref{Lvtau} fails for a given $v$ is less than
\[
\Cc{|L_v''|}{J} (2D^{-\vs'})^J=o( n^{-1}).
\]
Here note $D^{-\vs' J} = n^{-2}$ while 
the binomial coefficient is less than about $(\vs'\log D)^J$.)]

Thus for the lists
$L_v''':=L_v''\sm \gs_B(B\cap N_v)$ ($v\in V''$), we have w.h.p.\ (namely, whenever 
\eqref{Lvtau} holds)
\beq{Lv'''}
|L_v'''|>(1-e^{-1} + \vt'-\vs-o(1))t \,\,\,\,\forall v\in V''.
\enq

\mn

We then need an adjustment to get rid of large color degrees.
Let
\beq{dgcv}
d_{\gc}(v) =|\{w\in N_v\sm (T\cup B):\gc\in L_w'\}|
\enq
(to be clear: we could replace
$L_w'$ by $L_w'''$ here)
and 
\beq{Qgcv}
Q_\gc(v) =\{d_\gc(v)> (1-e^{-1} +\vs)t\}.
\enq
Theorem~\ref{T2.1} (with \eqref{ideal1}) gives (again with $\vs'\approx \vs^2$)
\beq{prQ}
\pr(Q_\gc(v)) < n^{-\vs'},
\enq
which, we assert, implies that w.h.p.\ 
\beq{badgc}
|\{\gc\in L_v': \mbox{$Q_\gc(v)$ holds}\}| < 2/\vs' \,\,\,\, \forall v \in V''.
\enq

\nin
\emph{Proof.}  For a given $v$, \eqref{prQ} implies
\[
\pr(|\{\gc\in L_v': \mbox{$Q_\gc(v)$ holds}\}| \geq 2/\vs')<
\C{|L_v'|}{2/\vs'}n^{-\vs'\cdot(2/\vs')} =o( n^{-1}).
\]
The first (intuitively clear) inequality follows from
\emph{negative association} (NA) of the indicators
$X_{w,\gc}:=\mathbbm{1}_{\{\gc\in L_w'\}}$ ($w\in V'$, $\gc\in \gG$).
See e.g.\ \cite{Pemantle} for the definition.  
Here NA follows from two easy observations:
(i) for a fixed $w$ the $X_{w,\gc}$'s are negatively associated 
(see \cite{BJ} for a more interesting statement); 
(ii) if each of the independent collections $\{X_{i,j}\}_j$ is negatively 
associated, then $\{X_{i,j}\}_{i,j}$ 
is as well; see e.g.\ \cite[Prop.\ 8]{DR}.

\qed

Finally, if \eqref{Lv'''} and \eqref{badgc} hold, then the lists 
\[
L_v'''\sm \{\gc\in L_v': \mbox{$Q_\gc(v)$ holds}\}  \,\,\,\,(v\in V'') 
\]
have sizes as in \eqref{Lv'''} (with the lost $2/\vs'$ from \eqref{badgc}
absorbed in the $o(t)$)
and color degrees at most $(1-e^{-1} +\vs)t$; so 
(since $\vs =\vt'/3$)
Theorem~\ref{TRS} gives the desired coloring of $V''$.

\section{Martingales}\label{Mart}

As mentioned earlier, Lemmas~\ref{LPh1a} and \ref{LC0} are 
based on martingale concentration results;  here we give
a reasonably convenient framework, mostly extracted from \cite[Sec.\ 3]{aglc},
that suffices for present purposes.

We are interested in the concentration of a r.v.\ $X =X(\go)$, with
$\go=(\go_1\dots \go_m)$ drawn from some
$\gO =\prod_{i\in [m]}\gO_i$.
We will need two versions, the first of which assumes a product measure on $\gO$:

\begin{itemize}

\item
$\,\,\,\go_1\dots\go_m$ are independent, with 
$\go_i$ uniform from $\gO_i$;

\item
for each $i\in [m]$ and $\eps=(\eps_1\dots \eps_m)\in \gO$ we are given some
$\gS_i=\gS_i(\eps_1\dots\eps_{i-1})\sub \gO_i$.

\end{itemize}

\nin
The second version---for use when we come to choosing $L_w^0$'s
(in Section~\ref{PLPh1a})---assumes:

\begin{itemize}

\item
$[m]:=S\times [\ell]$ ($S$ some set), ordered lexicographically;

\item
each $\go_{z,k}$ is uniform from $\gO_z\sm \{\go_{z,1}\dots \go_{z,k-1}\}$,
and is independent of all $\go_{x,j}$ with $x\neq z$;

\item
$X(\go)$ depends only on the \emph{sets} $\{\go_{z,1}\dots \go_{z,\ell}\}$;

\item
for each $z\in S$ we are given $\gS_z\sub \gO_z$ and for $i=(z,k)$
set $\gS_i =\gS_z$.

\end{itemize}

\nin
We will refer to these two situations as (A) and (B) respectively.
The following lemma applies to both.

\begin{lemma}\label{Laglc}
If 
\beq{EXX'}
|X(\eps_1\dots \eps_{i-1},\eps_i^0,\eps_{i+1}\dots \eps_m)-
X(\eps_1\dots \eps_{i-1},\eps_i^1,\eps_{i+1}\dots \eps_m)| 
\,\,\,
\left\{\begin{array}{ll}
\leq c,\\
=0 &\mbox{if $\eps_i^0,\eps_i^1\not\in \gS_i$}
\end{array}\right.
\enq
whenever $\eps_i^0,\eps_i^1\in \gO_i$ and $\eps_j\in \gO_j$ for $j\neq i$, and
\beq{PXEX}
\sum_i\pr(\go_i\in  \gS_i|\mbox{$\go_j=\eps_j$ for $j<i$}) \leq M,
\enq

\nin
then for any $\gl \in [0,16cM]$,
\beq{martbd}
\max\{\pr(X-\E X > \gl), \pr(X-\E X <- \gl)\}
 < \exp\left[-\frac{\gl^2}{32c^2 M}\right].
\enq
\end{lemma}
\nin
(The restriction on $\gl$ will never be a problem, and 
will usually be left to the reader.)

\mn
\emph{Note on} (A).
The conditioning in \eqref{PXEX} looks superfluous since the $\go_i$'s
are independent, but is needed since $\gS_i$ is a function of $\go_1\dots \go_{i-1}$.

\nin
\emph{Notes on} (B).  Here we will always have $S\sub V$, $\gO_z=\gG$ $\forall z$,
and $\{\go_{z,1}\dots \go_{z,\ell}\}=L_z^0 $
(so we think of choosing $L_z^0$ sequentially), and may use Lemma~\ref{Laglc} with
\beq{M}
M =\ell\sum_{z\in S} |\gS_z|/(D-\ell),
\enq
or any convenient upper bound on this
(the summand in \eqref{PXEX} for $i=(z,k)$ being at most
$|\gS_z|/(D-k+1)$).
In practice we define $X$ 
only on strings that can actually occur (those with $\go_{z,1}\dots \go_{z,\ell}$ 
distinct for each $z$).

\mn
\emph{Perspective.}  We think of the denominator in \eqref{martbd} as a variance
(of the sum of the martingale differences $X_i-X_{i-1}$, where $X_i=\E[X|\go_1\dots \go_i]$,
though here we refer to \cite{aglc}---also, in part, to the passage from \eqref{EXX'} to
\eqref{EXX''} shown below---for the connection); so \eqref{martbd} looks like the Central Limit Theorem.

Regarding \eqref{EXX'} and the meaning of the $\gS_i$'s, we can manage only a hint:
The first bound here is familiar from ``Azuma-type'' concentration bounds (e.g.\ 
\cite{McD}), and bounds the ``variance'' by $mc^2$.  This is natural for a sum
of independent r.v.s $\xi_i$ with absolute values at most $c$, but is weak if $\xi_i$ is mostly 
concentrated---say with probability $1-p_i$---on some single value, in which case the variance 
bound drops to $\sum p_ic^2$.  This will often be closer to our situation, with the 
$\gS_i$'s corresponding to (again, we omit the connection) the sets of ``exceptional'' values.

(The more significant point from \cite{aglc}---see its Lemma~3.4---is that one can allow
the r.h.s.\ of \eqref{EXX''} to depend on the history $\eps_1\dots \eps_{i-1}$.)

\mn
\emph{Derivation.}
As already said, most of Lemma~\ref{Laglc} derives from \cite{aglc}, 
and, not wanting to repeat that discussion here, we give only brief indications.

Lemmas~3.4 and 3.7 of \cite{aglc} (with the usual 
$\pr(Z>\gl) < e^{-\gz \gl}\E e^{\gz Z}$) imply
Lemma~\ref{Laglc} with \eqref{EXX'} replaced by its expectation version, namely
\beq{EXX''}
|\E [X|\eps_1\dots \eps_{i-1},\eps_i^0]-
\E [X|\eps_1\dots \eps_{i-1},\eps_i^1]| 
\,\,\,
\left\{\begin{array}{ll}
\leq c\\
=0 &\mbox{if $\eps_i^0,\eps_i^1\not\in \gS_i$;}
\end{array}\right.
\enq
so we need to say that \eqref{EXX''} follows from \eqref{EXX'}.
For (A) this is shown in \cite{aglc} (see Lemma~3.9(a)), but the argument is standard.
The corresponding argument for (B) is similar and presumably known, but is given here 
since we don't know a reference.

Fix $i=(z,k)$, let $\mu$ and $\nu$ range over $\prod_{j<i}\gO_j$ and $\prod_{j>i}\gO_j$
respectively, and  rewrite \eqref{EXX''} as
\beq{sumnu}
|\sum_\nu\left\{\pr(\nu|\mu\eps_i^0)X(\mu\eps_i^0\nu)
-\pr(\nu|\mu\eps_i^1)X(\mu\eps_i^1\nu)\right\}|
\,\,\,
\left\{\begin{array}{ll}
\leq c\\
=0 &\mbox{if $\eps_i^0,\eps_i^1\not\in \gS_i$ $~ (=\gS_i(\mu)$).}
\end{array}\right.
\enq
Set $\nu_z=(\nu_{z,k+1}\dots \nu_{z,\ell})$ and notice that

\nin
(i)  if $\eps_i^0$ appears in $\nu_z$ then
$\pr(\nu|\mu\eps_i^0) = 0$, and similarly with $\eps_i^0$
replaced by $\eps_i^1$;

\nin
(ii)
if $\eps_i^1$ appears in $\nu_z$ and $\eps_i^0$ does not, and 
$\nu'$ is $\nu$
with $\eps_i^1$ replaced by $\eps_i^0$ in $\nu_z$, then
(in view of the invariance of $X$ under permutations of $\{\go_{z,1}\dots\go_{z,\ell}\}$)
\[ 
\pr(\nu|\mu\eps_i^0)X(\mu\eps_i^0\nu)
=\pr(\nu'|\mu\eps_i^1)X(\mu\eps_i^1\nu').
\]

\nin
Thus the l.h.s.\ of \eqref{sumnu} doesn't change if we restrict the sums to $\nu$ with
neither of $\eps_i^0$, $\eps_i^1$ appearing in $\nu_z$; so \eqref{sumnu} follows from \eqref{EXX'}.\qed

\section{Concentration}\label{Conc}

The eventual goal here, addressed in Section~\ref{P23}, is establishing
\eqref{ideal1} and \eqref{ideal2}; but we first (in Section~\ref{PLPh1a})
dispose of Lemma~\ref{LPh1a}, which underlies the analysis for ``large'' $n$.

In Section~\ref{PLPh1a} and the ``small $n$'' part of Section~\ref{P23},
we aim for union bounds over sets of vertices of bounded size,, so want to say
that each of various desired events fails with probability $n^{-\go(1)}$;  in these cases we
use``w.v.h.p.'' to mean with probability $1-n^{-\go(1)}$.

\subsection{Proof of Lemma~\ref{LPh1a}}\label{PLPh1a}
This is all based on (regime (B) of) Lemma~\ref{Laglc}, though, to be clear,
the use of this machine for
\eqref{whp3} and \eqref{whp4} is pretty excessive.
In the proofs of \eqref{whp1} and \eqref{whp2} we use  
$N^0_\gc(v) = \{w\sim v:\gc\in L_w^0\}$;
$N(W) =\cup_{w\in W}N_w$ (and similarly for $N_\gc^0(W)$); and
$d_\gc^0(\cdot) =|N^0_\gc(\cdot)|$.

\mn
\emph{Proof of \eqref{whp1}.} Here \eqref{EEP} (with \eqref{przT'}) gives 
$\E'\pr''(z\in T)=\pr(z\in T)\sim e^{-1}$,
so we just need concentration of 
\beq{przT}
\pr''(z\in T) =\sum_{\gc\in L_z^0}\ell^{-1}(1-\ell^{-1})^{d_\gc^0(z)}.
\enq

Given $L_z^0$, the sum in \eqref{przT} is a function of the events
$\{\gc\in L_y^0\}$,
with 
$\gc$ running over $ L_z^0$ and $y$ over $N_z$.  
We 
apply Lemma~\ref{Laglc} (as sketched in the paragraph containing \eqref{M}),
with
\beq{firstuse}
\mbox{$S=N_z$,  $~\gO_y=\gG$,  $~\gS_y=L_z^0$, $~c=\ell^{-2}~$
and $M = D\ell^2/(D-\ell)$.}
\enq~
[The second line of \eqref{EXX'} is obvious, the rationale for $c$ 
should be clear from \eqref{przT}, and $M$ is as in \eqref{M}.]

Thus for any
$\gl =\go(\ell^{-1/2})$, Lemma~\ref{Laglc} gives 
\[
\pr'(|\pr''(z\in T)- \E' \pr''(z\in T)|>\gl|L_z^0) < e^{-\go(\ell)}=n^{-\go(1)},
\]
implying \eqref{whp1}.\qed

\mn
\emph{Proof of \eqref{whp2}.}
Let 
\beq{fW}
f^v_{wz}
 = f^v_{wz}
(L_{wz}^0) =\E' [\pr''(A^v_{wz})|L_{wz}^0]
\enq
(for distinct $w,z\in N_v$, where $L_{wz}^0=L_w^0\cap L_z^0$) and
\[
\mbox{$X_v = \sum\{ f^v_{wz}:\{w,z\}\sub N_v\}$}.
\]

Again using \eqref{EEP}, we have 
\beq{Exps}
\pr(A^v_{wz})=\E'\pr''(A^v_{wz}) = \E'\{\E'[\pr''(A^v_{wz})|L_{wz}^0]\} =\E' f^v_{wz}
\enq
and (with sums over $\{w,z\}\sub N_v$ with $w\not\sim v$)
\begin{eqnarray}
|\sum(\pr''(A^v_{wz}) -\pr(A^v_{wz}))| &=&
|\sum(\pr''(A^v_{wz}) -\E' f^v_{wz})| 
\nonumber\\
&=&
|\sum(\pr''(A^v_{wz}) -f^v_{wz} + f^v_{wz}- \E' f^v_{wz})| \nonumber\\
&\leq & 
\sum|\pr''(A^v_{wz}) -f^v_{wz}| + |X_v- \E' X_v| .\label{P''AWfW}
\end{eqnarray}

\nin
So we just need to show 
\beq{toshow}
\mbox{each of the two parts of \eqref{P''AWfW} is $o(D)$ w.v.h.p.}
\enq

For the first part of this, recall from \eqref{Jvwz}
that $J^v_{wz} = N(v,w,z)\sm \{w,z\}$ and set
\[
J^0_\gc(v,wz) =\{x\in J^v_{wz}: \gc\in L_x^0\}.
\]
It turns out---unnecessarily here, but at no extra cost---that each 
\beq{prAwz}
\pr''(A^v_{wz}) = \sum_{\gc\in L_{wz}^0} \ell^{-2}(1-\ell^{-1})^{|J^0_\gc(v,wz)|}
\enq
is well concentrated; so we work with these.

Given $L_{wz}^0$,
we apply Lemma~\ref{Laglc} to $X:=\pr''(A^v_{wz})$
with
\beq{seconduse}
\mbox{$S=J^v_{wz}$,  $~\gO_y=\gG$,  
$~\gS_y=L_{wz}^0$, $~c=\ell^{-3}~$ and $M = 4\ell^2$.}
\enq
[The validity of $c$ should again be clear (from \eqref{prAwz}), and 
$M\geq \ell |J^v_{wz}||L_{wz}^0|/(D-\ell)$
(the value in \eqref{M}).]

Thus, for any given $L_{wz}^0$ (and $\gl\in [0,64\ell^{-1}]$),
\[
\pr'(|\pr''(A^v_{wz})-f^v_{wz}| > \gl|L_{wz}^0) < \exp[-\gO(\gl^2 \ell^4)],
\]
which is $n^{-\go(1)}$ if $\gl =\go(\ell^{-3/2})$.
So if $\gl$ is both $\go(\ell^{-3/2})$ and $o(\ell^{-10/9})$ (see \eqref{nlarge}),
then
w.v.h.p.\
\[
|\pr''(A^v_{wz})-f^v_{wz}| < \gl\,\,\,\,\forall w,z
\] 
and 
\[
\sum_{\{w,z\}\sub N_v}|\pr''(A^v_{wz})-f^v_{wz}| < \gl \C{D}{2} =o( D),
\]
yielding the first part of \eqref{toshow}.

\mn

We turn to the second part of \eqref{toshow}, observing for perspective that
here we \emph{don't} expect good concentration for individual $f^v_{wz}$'s 
(since $f^v_{wz} =\Theta( |L_{wz}^0|\ell^{-2})$ and the $|L_{wz}^0|$'s are
not very concentrated).

Noting that $X_v $ is a function of $L^0|_{N_v}$,
we apply Lemma~\ref{Laglc}
with
\[ 
\mbox{$S=N_v$,  $~\gO_z=\gS_z=\gG$,  $~c=D\ell^{-2}~$
and $~M = D\ell$.}
\]  
[For $M$ we just bound the summands in \eqref{PXEX} by 1.
To see that \eqref{EXX'} holds, notice that, for $i=(z,j)$, the only possible
effects of changing $\go_i$ from $\ga$ to $\gb$
(while holding other entries fixed, as in \eqref{EXX'}) are removal of $\ga$ from and/or
addition of $\gb$ to some of the sets $L_{wz}^0$ (with $w\in  N_v$),
and each such modification changes the corresponding $f^v_{wz}$ by at most $\ell^{-2}$
(cf. \eqref{prAwz}).]

This gives 
\[
\pr'(|X_v-\E' X_v|>\gl ) < \exp[-\gO(\gl^2 (\ell/D)^3)] ,
\]
and using this with $\gl$ satisfying $\go(D^{3/2}\ell^{-1})=\gl =o(D)$
(possible since we assume \eqref{nlarge}) 
gives the second part of \eqref{toshow}.

\qed

\nin
\emph{Proof of \eqref{whp3}.}
Fix $v,Y, \gG_0$ as in \eqref{whp3}. For the quantity of interest,
$X:=\sum_{y\in Y} |L_y^0\cap \gG_0|$, we have
\[
\E' X =\theta D\cdot \theta D/C\cdot \ell/D = \theta^2D\ell/C,
\]
so again just need concentration. 
We use Lemma~\ref{Laglc} with 
\[    
\mbox{$S=Y$,  $~\gO_z=\gG$,  $~\gS_y=\gG_0$, $~c=1~$
and $~M = |Y|\ell\cdot |\gG_0|/(D-\ell)< (1+\gd) \theta^2D\ell/C$,}
\]    
yielding (with plenty of room)
\beq{X2X}
\pr'(X> 2\E' X) < \exp[-\gO(\theta^2 D\ell/C)],
\enq
which, since $C\theta^{-2}=o(\ell)$ (see following \eqref{Ceps})
is still $n^{-\go(1)}$ when summed over the fewer than 
$n2^D$ choices for $(v,Y,\gG_0)$.

\qed

\mn
\emph{Proof of \eqref{whp4}.}
Fix $(v,\tau)$ and let 
$X = \sum_{w\not\sim v}|\tau(N_w\cap N_v)\cap L_w^0|$
(the sum in \eqref{whp4}), noting that
\[
\E' X ~=~ \sum_{w\not\sim v}|\tau(N_w\cap N_v)|\ell/D ~\leq ~
 \sum_{w\not\sim v}|N_w\cap N_v|\ell/D.
 \]
(Of course $w$'s at distance more than 2 from $v$ play no role here.)
For Lemma~\ref{Laglc} we take
\[
\mbox{$S=\{w:\textrm{dist}(v,w)=2\}$,  $~\gO_z=\gG$,  
$~\gS_w=\tau(N_w\cap N_v)$, $~c=1$, }
\]
and 
\[
M = D^2\ell/(D-\ell) \geq \ell\sum_{w\not\sim v}|N_w\cap N_v|/(D-\ell)
\geq \ell\sum_{w\not\sim v}|\tau(N_w\cap N_v)|/(D-\ell)
\]    
(the r.h.s.\ being the value in \eqref{M}).  This gives
\[
\pr'(X>2D\ell) < \exp[-\gO(D\ell)],
\]
and \eqref{whp4} follows (the number of $(v,\tau)$'s, roughly $D^D$,
again being negligible by comparison).

\qed

\subsection{Properties \eqref{ideal1} and \eqref{ideal2}}\label{P23}  
Here we prove \eqref{small23} and \eqref{CL0}; so (e.g.) ``proof of 
\eqref{small23}'' is shorthand for ``proof that \eqref{small23} holds with
[the appropriate probability].''
Here we find it convenient to take $V=[n] $
and \emph{sometimes} use $i,j,\ldots $ instead of $v,w,\ldots$
for vertices.

We treat small and large $n$ separately, beginning with the former.
When we do come to large $n$, we are given $L_v^0$'s, which we may (and do) assume 
satisfy \eqref{whp1}-\eqref{whp4}.

For $\tau:N_v\ra\gG$, let
\beq{Z}
Z(v,\tau) =\{z\in N_v: |\tau^{-1}(\tau_z)\cap N_v|< C\} ,
\enq
and say $(v,\tau)$ is \emph{bad} if $|Z(v,\tau)| < (1-\theta)D$ (and \emph{good} otherwise).

Before attacking \eqref{ideal1} and \eqref{ideal2}, we show bad $(v,\tau_{N_v})$'s 
(for our random $\tau$) are unlikely;
this will be used to control the differences in \eqref{EXX'} of Lemma~\ref{Laglc}.
Fix $v$ and notice that if $(v,\tau|_{N_v})$ is bad
then there are
$Y\in \C{N_v}{\theta D}$ and $\gG_0\in \C{\gG}{\theta D/C}$ with $\tau(Y)\sub \gG_0$;
thus for small $n$,
\beq{Zbad}
\,\,
\pr(\mbox{$(v,\tau|_{N_v})$ bad}) < \C{D}{\theta D}\C{D}{\theta D/C}\left(\theta/C\right)^{\theta D}
< \left[(e/\theta)(eC/\theta)^{1/C}(\theta/C)\right]^{\theta D}<
e^{o(D)}(e/C)^{\theta D} \,\, = e^{-\go(D)}
\enq
(using \eqref{Ceps}), while \eqref{whp3} (with the AM-GM Inequality)
gives the corresponding bound for large $n$:
\begin{eqnarray}
\label{Zbad'}
\,\,\,\,\,
\pr''(\mbox{$(v,\tau|_{N_v})$ bad}) &<&
\sum_{\gG_0}\sum_Y\prod_{y\in Y}\frac{|L_y^0\cap \gG_0|}{\ell}\\
&\leq& \C{D}{\theta D/C}\C{D}{\theta D}\left(\frac{1}{\theta D}
\sum_{y\in Y} \frac{|L_y^0\cap \gG_0|}{\ell}\right)^{\theta D}
<\left(\frac{eC}{\theta}\right)^{\theta D/C}
\left(\frac{2e}{C}\right)^{\theta D} = e^{-\go(D)}.\nonumber
\end{eqnarray}

\mn

For the rest of this section
we fix $v\in V^*$ 
(as mentioned earlier, this restriction is really only needed for \eqref{ideal2})
and 
let $\tau= (\tau_1\dots \tau_n)$, with the $\tau_i$'s independent and
\[
\mbox{$\tau_i$ uniform from}\left\{\begin{array}{ll} 
\gG&\mbox{if $n$ is small,}\\
L_i^0&\mbox{if $n$ is large.}
\end{array}\right.
\]
For convenience we assume $N_v=[D]$ and in what follows
will use whichever of $N_v$, $[D]$ seems more appropriate;
we will always use $\tau_{[D]}$ for the restriction of $\tau$ to $N_v$.

The proofs of \eqref{ideal1} and \eqref{ideal2} are again applications of Lemma~\ref{Laglc}.
Here we are in regime (A) of the lemma and will always take $\gO=\prod_{i\in [n]}\gO_i$,
noting that most $i$'s---those with $\textrm{dist}(v,i)>2$---will never affect the $X$ in
question; such $i$'s will have $\gS_i=\0$, so won't affect $M$ (see \eqref{PXEX}) or 
the bound \eqref{martbd}.  (Alternatively, we could take $\gO=\prod_{i\in [m]}\gO_i$,
with $[m] =\{i:\textrm{dist}(v,i)\leq 2\}$.)

\nin
\emph{Proof of \eqref{ideal1} for small n.}
Let 
\[
X=X(\tau)=|T\cap Z|=\sum_{i\in N_v}X_i,
\] 
where $Z=Z(v,\tau_{[D]})$ and $X_i =\mbone_{\{i\in T\cap Z\}}$.

We first observe that 
\beq{EX}
\E X\sim e^{-1}D
\enq
and that we will have \eqref{ideal1}, i.e. 
\[
\mbox{ w.v.h.p. $~|T\cap N_v|\sim e^{-1}D~$},
\]
if we can show 
\beq{Xvtoshow}
\mbox{w.v.h.p.  $~|X-\E X| =o( D)$.}
\enq
[Because:
With $U:=N_v\sm Z$, the combination of \eqref{Zbad} and \eqref{Ceps} 
(with $\E |U| < \theta D + D\pr(|U|> \theta D)$)
gives 
\beq{EU}
\mbox{$\E |U| \less \theta D~$ and w.v.h.p. $|U| \less \theta D$,}
\enq
which, since ($0\leq $)   $|T\cap N_v|-X \leq |U|$, 
 implies \eqref{EX} (see \eqref{ideal1'}) and
\beq{TNv}
\mbox{w.v.h.p. $||T\cap N_v|-X| \less \theta D$.]}
\enq

\nin
\emph{Proof of \eqref{Xvtoshow}.}
Notice (with justification to follow) that, if
$\tau'_i=\gb$ and $\tau_j=\tau'_j$ $~\forall j\neq i$, then
\beq{Obstaus}
X(\tau)-X(\tau')~\leq~ \left\{\begin{array}{ll}
 C &\mbox{(always),}\\
0& \mbox{if $i>D$ and $i\not\sim \tau^{-1}(\gb)\cap Z$.}
\end{array}\right.
\enq
(Note $Z$ is determined by $\tau_{[D]}$, so its appearance here makes sense.)
This is true because, with $X=X(\tau)$, $X'=X(\tau')$, and $j\in [D]$,
we can only have $X_j> X'_j$ (i.e.\ $X'_j=0$, $X_j=1$)
if 

\begin{itemize}

\item
$i=j$, or

\item
$i\sim j$, 
$\tau_j=\gb$, and $|\tau^{-1}(\gb)\cap (N_v\sm\{i\})| <C$, or

\item
$i\leq D$, $\tau_j=\gb$,
and $|\tau^{-1}(\gb)\cap (N_v\sm\{i\})|=C-1$.

\mn
\end{itemize}

For Lemma~\ref{Laglc} we take

\begin{itemize}

\item
$\gO_i=\gG$, $\go_i=\tau_i$  ($\forall i$);

\item
for $i\in [D]$, $\gS_i=\gO_i$ (so the summand in \eqref{PXEX} is 1); 

\item
for $i>D$, $\gS_i= \gS_i(\tau_{[D]}) =\tau(N_i\cap N_v)$; and

\item
$c=C$.

\end{itemize}
Then \eqref{EXX'} is given by \eqref{Obstaus}, and
for $i>D$,
\beq{prgoi}
\pr(\go_i\in \gS_i)= |\tau(N_i\cap Z)|/D\leq |N_i\cap Z|/D,
\enq
implying that the sum in \eqref{PXEX} is at most 
\[
D + D^{-1}\sum_{i>D}|N_i\cap N_v| < 2D =:M.
\]

Thus, finally, Lemma~\ref{Laglc} gives
\[
\pr(|X-\E X|>\gl) =\exp[-\gO(\gl^2/(C^2D))];
\]

\nin
so $|X-\E X| < \gl$ w.v.h.p.\ for any $\gl=\go( C(D\log n)^{1/2})$,
and in particular, \emph{since $n$ is small}, for some $\gl =o(D)$.

\qed

\mn
\emph{Proof of \eqref{ideal1} for large n.}
Recall that the $\tau_i$'s are now drawn uniformly from the $L_i^0$'s, which 
are assumed to satisfy \eqref{whp1}-\eqref{whp4}, and that we aim for \eqref{CL0};
that is, we just need to say, using the Local Lemma, that with \emph{positive} probability
\eqref{ideal1} 
holds for all $v\in V^*$.
(Again, the restriction to $V^*$ is unnecessary here.)

Thus, since each event in \eqref{ideal1} 
is independent of all but $O(D^4)$ others
(since each is determined by the restriction of $\tau$ to vertices within distance 2 of $v$),
even failure probabilities below $D^{-5}$ will suffice here (though we will do much better); 
so we now use \emph{w.f.h.p.}\ (f = fairly) for ``with probability at least $1-D^{-5}$.''

The arguments here are similar to those for small $n$, so we mainly indicate differences.
As before, we have $v\in V=[n]$,
$N_v=[D]$ and $\tau =(\tau_1\dots \tau_n)$, 
and set 
\[
X=X(\tau)=|T\cap Z|=\sum_{z\in N_v}X_z,
\] 
where $Z=Z(v,\tau_{[D]})$ and $X_z =\mbone_{\{z\in T\cap Z\}}$.

As was earlier true of \eqref{Xvtoshow}, in this case it is enough to show
\beq{Xvtoshow'}
\mbox{w.f.h.p.\  $|X-\E'' X| =o( D)$.}
\enq
[The justification for ``enough to show'' follows the earlier one:  we now add \eqref{EEP} to \eqref{ideal1'}
to give \eqref{EX}, and \eqref{EU} is given by \eqref{Zbad'} instead of \eqref{Zbad}.
(``W.v.h.p.'' is still okay in both of these, but not in \eqref{Xvtoshow'} and \eqref{ideal1}.)]

We again have \eqref{Obstaus} (with the same justification); to repeat:
\beq{Obstaus'}
X(\tau)-X(\tau')~ \left\{\begin{array}{ll}
\leq C \\
\leq 0& \mbox{if $i>D$ and $i\not\sim \tau^{-1}(\gb)\cap Z$}
\end{array}\right.
\enq
whenever $\tau'_i=\gb$ and $\tau_j=\tau'_j$ $~\forall j\neq i$.

The parameters for Lemma~\ref{Laglc} are also as before, except
$\gO_i$ changes from $\gG$ to $L_i^0$; thus:

\begin{itemize}

\item
$\gO_i=L_i^0$, $\go_i=\tau_i$  ($\forall i$);

\item
for $i\in [D]$, $\gS_i=\gO_i$; 

\item
for $i>D$, $\gS_i= \gS_i(\tau_{[D]}) =\tau(N_i\cap N_v)$; and

\item
$c=C$.

\end{itemize}

The biggest (if still minor) change here is that the 
$|\tau(N_i\cap Z)|/D$
of \eqref{prgoi} becomes 
$|\tau(N_i\cap Z)\cap L_i^0|/\ell$;
that is, for $i>D$,
\[
\pr''(\go_i\in \gS_i(\tau_{[D]}))= |\tau(N_i\cap Z)\cap L_i^0|/\ell,
\]
so the sum in \eqref{PXEX} is at most 
\[
D + \ell^{-1}\sum_{i>D}|\tau(N_i\cap Z)\cap L_i^0| < 3D =:M,
\]
with the inequality given by \eqref{whp4}.
So Lemma~\ref{Laglc} again gives
\[
\pr''(|X-\E X|>\gl) =\exp[-\gO(\gl^2/(C^2D))],
\]
which in particular says 
$|X-\E X| =o( D)$ w.f.h.p. 
(The earlier ``w.v.h.p.'' needed $\gl$ with
$\go(C\sqrt{D\log n})=\gl=o( D)$, which may not exist when $n$ is large.)

\qed

\mn
\emph{Proof of \eqref{ideal2'} for small and large n.} 
Recall this refers to the parts of \eqref{small23} and \eqref{CL0} concerning \eqref{ideal2'}.
There is nothing new here, so we keep it brief.

Fix $v$ and let
\[
X=|\{wz\in \ov{G}[N_v]:\mbox{$A_{wz}^v $ holds}\}|.
\]
By \eqref{Ewz} if $n$ is small, or \eqref{Ewz} and \eqref{EEP} if $n$ is large,
\[
\E X > \vt De^{-3},
\]
so it is enough to show
\beq{XEXD}
|X-\E X| =o( D) \,\,\left\{\begin{array}{ll}
\mbox{w.v.h.p.} & \mbox{if $n$ is small,}\\
\mbox{w.p.p.} & \mbox{if $n$ is large.}
\end{array}\right.
\enq

The argument(s) for this are nearly identical to those for \eqref{ideal1} above, 
the sole diffence being that $C$ (so also $c$) is replaced by $1$.
(In words:  changing a single $\tau_i$ can't change the present $X$ by more than 1.)
Thus Lemma~\ref{Laglc} gives (slightly improving the earlier bounds)
\[
\pr(|X-\E X|>\gl) =\exp[-\gO(\gl^2/D)],
\]
and \eqref{XEXD} follows.

\qed

\section{Clusters}\label{Clusters}

We now assume $V^*$ has been colored by $\gs$, and choose the ``real'' lists $L_v$ 
for $v\in \cup C_i$.
(Thus, as promised, we discard the earlier ``dummy'' colors.
We again note that the following analysis is valid for \emph{any} $\gs:V^*\ra\gG$.
As noted at the end of Section~\ref{Sketch}, we now \emph{occasionally}---only where 
relevant---use $\DD$ for $|\gG|$ ($=D+1$) in place of the slightly incorrect $D$.)

With $V_i =V^*\cup \cup_{j\le i}C_j$, 
it is enough to show that for $i=1\ldots $ and any coloring $\gs$ of $V_{i-1}$
(note we continue to use $\gs$ for the growing coloring),
\beq{notext}
\pr(\mbox{$\gs$ does not extend to $C_i$}) =o( D/n).
\enq
We prove this here modulo a few assertions 
that will be justified in
Sections~\ref{PLsucceed} and \ref{Matchings}.

\nin
\emph{Note on parameters.}
For the rest of our discussion, $\eps$, though fixed, may best be thought of as tending slowly to 0.  
(Letting it \emph{actually} tend to 0 would be more convenient here but complicates the 
presentation of the sparse phase.)  We will use $a\ll b$ as a finitary counterpart of $a=o(b)$, meaning:  
for small enough $\eps$ we can choose the other parameters below so that any use of 
$a\ll b$ promises that (say, though this is more than we need) $a< 0.01 b$.  Of course this will sometimes
imply larger gaps; e.g.\ existence of $\gz_0$ in \eqref{epszeta} requires $\eps < 0.0001$.
We then replace  $\sim$
by $\approx$ (and $\less$ by $\lesser$); e.g.\
$a\approx b$ means $|a-b|\ll a$.
(Recall from ``Usage'' in Section~\ref{Sketch} that all our assertions are for sufficiently large $D$.)
Finally, we recall (from near the end of Section~\ref{Sketch}) that the implied constants in
$O(\cdot)$, $\gO(\cdot)$, $\Theta(\cdot)$
do not depend on our various parameters.

\mn
\emph{Proof of \eqref{notext}.}
Write $C$ for $C_i$, recalling that
(for convenience slightly relaxing from Lemma~\ref{LACK})
\beq{EDepsD}
\mbox{$||C| -D| < b\eps D;~$ $~C$ has external degrees and
internal \emph{non}-degrees at most $b\eps  D$,}
\enq  
with $b$ fixed (e.g.\ $b=7$).
Say $C$ is \emph{large} if $|C|> |\gG|$ ($=\DD$)
and \emph{small} otherwise.

Let $\gz_0$ satisfy 
\beq{epszeta}
\eps \ll \gz_0 D \ll 1.
\enq
(We could, e.g., specify $\gz_0 = \sqrt{\eps}/D$, but \eqref{epszeta} is all we'll use.)
Let $H=\ov{G}[C]$ and $|H| =\gz D^2$,
and 
say $\gz$ is \emph{small} if $\gz < \gz_0$ and \emph{large} otherwise.

Notice that $\gz $ small implies $C$ is small (since $C$ large implies $\gz > 1/(2D)$).
In this case we set $C'=C$ and $\gG'=\gG$, and jump to the paragraph 
containing \eqref{R'}.  If $\gz$ is large, we first perform the following pairing step
(at ``Process''),
which, as suggested in the ``dense phase'' preview of
Section~\ref{Sketch}, will provide a little room
in the coming use of Hall's Theorem by arranging repeated use of a few colors.
This requires a few parameters.

Let $\eta =\eta_\gz $ satisfy
\beq{etazeta}
\max\{\gz,1/D\} \ll \eta \ll \gz/\eps \,\,\,\,(=O(1))
\enq
and set
\beq{q}
q ~= ~1-\exp\left[-\frac{\gz\kkk}{8b\eps}\right] ~\gg~\eta,
\enq
where $b$ is from \eqref{EDepsD} and
the inequality uses $\gz\kkk /\eps\gg \eta\kkk$ 
(see \eqref{etazeta})
if $\gz \kkk <0.1b\eps$ (say), and $\eta \ll 1$ otherwise, and
\beq{KetaDq}
m = K\eta D/q,
\enq
with $K$ satisfying
\beq{K}
\eta \ll K\eta \ll q\ll (\gz D/\eps) K\eta
\enq
(possible since $q\gg \eta$ and $\gz\gg\eps/D$
[see \eqref{q} and \eqref{epszeta})]),
so in particular 
\beq{mllD}
m\ll D.
\enq

Set 
\beq{fgc}
f(\gc)=
\sum\{\gz D +d_H(v):v\in C, v\sim\gs^{-1}(\gc)\}
\enq
and let $\gG=\{\gc_1,\dots\gc_\DD\}$, with $f(\gc_1)\leq f(\gc_2)\leq \cdots$.
\begin{obs}\label{Obsf}
For any $i\leq m$, 
$
\,\,f(\gc_i) =O(\eps \gz D^2).
$
\end{obs}

\nin
\emph{Proof.}
It is enough to show
\beq{fsum}
\sum_\gc f(\gc) =O(\eps \gz D^3),
\enq
since this, with \eqref{mllD}, gives 
$
f(\gc_i)\leq f(\gc_m) < (D-m)^{-1} \sum_\gc f(\gc) =O(\eps \gz D^2).
$

For \eqref{fsum}, we just observe that the contribution of the $\gz D$ terms (in \eqref{fgc}) 
to the sum in \eqref{fsum} is at most 
\[
\gz D |\nabla (C)| =\gz D\cdot O(\eps D^2),
\]
while that of the $d_H(v)$ terms is 
\[
\sum_{v\in C}d_H(v)|\{\gc: v\sim \gs^{-1}(\gc)\}| \leq b\eps D\cdot 2\gz D^2
\]
(since $\sum d_H(v) = 2\gz D^2$ and \eqref{EDepsD} says
$|\{\gc: v\sim \gs^{-1}(\gc)\}| \leq b\eps D$).
\qed

\mn

Noting that clusters other than $C$ play no further role here, we now repurpose the notation $C_j$:

\mn
\textbf{Process}: 
Set $C_0=C$, $H_0=H$ and $ M_0=\0$, and
for $i=1\dots m$, 
reveal 
\beq{Ji}
J_i:=\{v\in C_{i-1}\sm N(\gs^{-1}(\gc_i)):\gc_i\in L_v\}
\enq 
and do:

\begin{itemize}

\item[(I)]
If $H[J_i]\neq \0$, choose (arbitrarily) $x_iy_i\in H[J_i]$ and set
\[ 
\gs(x_i)=\gs(y_i)=\gc_i,  \,\, M_i=M_{i-1}\cup \{x_iy_i\}, \,\,C_i=C_{i-1}\sm \{x_i,y_i\},
\,\,\text{and} \,\, H_i =H_{i-1}-\{x_i,y_i\} \, \,(=H[C_i]).
\]

\item[(II)]
If $H[J_i]= \0$ and $J_i\neq \0$, choose $z_i\in J_i$ with $d_H(z_i)$ minimum and set
\[ 
\gs(z_i)=\gc_i,  \,\, M_i=M_{i-1}, \,\,C_i=C_{i-1}\sm \{z_i\},
\,\,\text{and} \,\, H_i =H_{i-1}-z_i \, \,(=H[C_i]).
\]

\item[(III)]
If $J_i=\0$ do nothing. (Updates are unnecessary here as we will take the process to have failed.)

\end{itemize}

\mn

Say the process \emph{succeeds} if

\nin
(S1)\, $|L_v\cap \{\gc_1\dots \gc_{m-1}\}| \leq 0.1 \gd \log n \,\,\, \forall v\in C$,

\nin
(S2)\, $J_i\neq \0$ $\,\forall i\in [m]$, and

\nin
(S3)\, $|M_m|\geq \eta D$.

\begin{lemma}\label{Lsucceed}
The process succeeds with probability $1-o(D/n)$.
\end{lemma}
\nin
\emph{Perspective.}
1.  If the process does succeed, then 
$|\gG|$
shrinks by $m$ and $|C|$ by the larger $m+|M_m|\geq m+\eta D$,
giving the promised room in Hall.
The natural requirement (in \eqref{mllD}) that $m$ be small relative to $D$ ensures that
the process removes only a negligible portion of $\gG$.

\nin
2.
The division into large and small $\gz$ turns out to be fairly delicate; roughly:
(a) there are situations with $\gz=O(\gd/D)$ where we cannot ``promise'' a full matching of
$C$ and $\gG$, so need the process to provide at least some repetition of colors as in (I);
while (b) if $\gz =\eps/D$, we \emph{cannot} promise that the process succeeds
(but eventually \emph{can} promise a full matching).
To see these, let $|C|=\DD$ ($=|\gG|$) and (with $v$ and $\gc$ always in $C$ and $\gG$ resp.)
consider:

\nin
(a) Suppose $\gz =\gd/D$ and let $H$ be a matching of size $\gz D^2=\gd D$.  
Then we can have $\gc$ with 
$|\{v:v\sim \gs^{-1}(\gc)\}| =|\{v:d_H(v)>0\}|=2\gd D$, and
\[
\pr(\gc\not\in L_v~\forall v\not\sim \gs^{-1}(\gc)) \approx (1-k/D)^{(1-2\gd )D} \approx n^{-(1-\gd)},
\]
which is not enough for \eqref{notext} unless $D$ is quite large.  (More trivially, $\gz =\DD/D^2$ allows
$\gc$ with $v\sim \gs^{-1}(\gc) $ $\forall v$ and \emph{no} possibility of the desired matching.)

\nin
(b)  Suppose $\gz=O(\eps /D)$ and let $H$ be the disjoint union of $\gz D/\eps=O(1)$ 
stars of size $\eps D$, with $T$ the set of centers.
Then the probability that the
process fails---even in the strong sense that $H[J_i]=\0$ $\forall i\in [m]$---is at least 
$\pr(J_i\cap T=\0 ~\forall i\in [m]) = (1-O(k/D))^m = n^{-O(m/D)}$, which, since $m\ll D$, 
is too large for Lemma~\ref{Lsucceed} (even with $o(1/n)$ replaced by $o(D/n)$ if $D$
isn't too large).

The proof of Lemma~\ref{Lsucceed} is given 
in Section~\ref{PLsucceed};
for now we assume that either $\gz$ is small---so we skipped ``process''---or  
the process has finished successfully (producing an updated $\gs$), and continue.

\mn

Write $C'$ and $\gG'$ for what's left of $C$ and $\gG$.
(So $C'=C$ and $\gG'=\gG$ if $\gz$ is small,
and $C'=C_m$, $\gG'=\{\gc_{m+1}\dots \gc_{\DD}\}$ if $\gz$ is large.)
Let $|C'| = J'$ and $|\gG'|=J'+R'$, noting that 
\beq{J'D}
J'\approx D
\enq 
(see \eqref{mllD}) and,
if $|C|=\DD \pm \ell$ with $\ell \geq 0$ (so $|C|=\DD -\ell$ iff $C$ is small),
\beq{R'}
R'\left\{\begin{array}{ll}
=\ell &\mbox{if $\gz$ is small (so $C$ is small),}\\
\ge \ell +\eta D&\mbox{if $\gz$ is large and $C$ is small,}\\
\ge \eta D-\ell \approx \eta D&\mbox{if $C$ is large.}
\end{array}\right.
\enq
[For the final ``$\approx$'': $\,\ell D/2\leq |H|=\gz D^2\ll \eta D^2$ implies $\ell\ll \eta D$.]
In particular $R'\geq 0$.

\mn

Set (as before; see \eqref{tgd'}) $\gd' = 0.9\gd$ and $t=(1+\gd') \log n$.
We now choose lists $L_v'$ ($v\in C'$) uniformly and independently
from $\C{\gG'}{t}$ for $v\in C'$, and want to show
\beq{Lv'C'}
\mbox{with probability $1-o(D/n)$ the $L_v'$'s support an extension of $\gs$ to $C'$.}
\enq
(This is enough since (S1) implies $|L_v\cap \gG'|\geq t$ $\forall v\in C$.  
Of course when $\gz$ is small, the $t$ here could be $k$.)

\nin

Let $B'$ be the $(C',\gG')$-bigraph of allowed colors (i.e.\ with $v\gc\in B'$ iff 
no neighbor of $v$ is currently colored $\gc$) 
and $K'$ the random subgraph with $v\gc\in K'$ (for $v\gc\in B'$) iff
$\gc\in L_v'$. 
Then for \eqref{Lv'C'} it is enough to show
\beq{notext'}
\pr(\mbox{$K'$ has no $C'$-perfect matching})= o( D/n)
\enq
(where a matching of a bigraph on $X\cup Y$ is \emph{X-perfect} if it uses all vertices of $X$).
For the proof of this we need the following supporting facts, which will be established 
in Section~\ref{Matchings}.  
In these, $t$ is as above, as are $\eta=\eta_\gz$ (see \eqref{etazeta})
and small/large for $\gz$ (see following \eqref{epszeta}).

\mn
\begin{lemma}\label{LMg}
Assume $B$ is bipartite on $X\cup Y$ (and $\oB$ is its bipartite complement),
with $|X|=J\approx D$; $|Y|=J+R$; $d_{\oB}(v) =O(\eps D)$  $\forall v\in X$; 
$|\oB| \leq 2\gz D^2$; and 

\nin
{\rm (a)}  $R=0\,$ if $\gz$ is small,

\nin
{\rm (b)}   $D\gg R\morer \eta D\,$ 
if $\gz$ is large. 

\nin
For each $v\in X$ choose $M_v $ uniformly from $\C{Y}{t}$ (these choices independent),
and let $K$ be the random subgraph of B with $v\sim_K\gc$\, iff \,
($v\sim_B$) $\gc \in M_v$.  Then 
$K$ admits an $X$-perfect matching with probability $1-o(D/n)$.
\end{lemma}
\nin
(The upper bound in (b) is unnecessary, but will make our arithmetical lives a little easier.)

Perhaps the most interesting point here is the next observation, which we will use
in handling small $C$.

\begin{lemma}\label{Lswitch}  
Let $F$ be bipartite on $X\cup Y$ with $|X|\leq |Y|$,
let $u,v\in Y$, and let $F'$ be obtained from $F$ by revising:  $N_{F'}(u) =N_{F}(u)\cup N_{F}(v)$
and $N_{F'}(v) =N_{F}(u)\cap N_{F}(v)$.
Let $L$ be the random bigraph on $X\cup Y$ gotten  by choosing neighborhoods
$L_x:=N_L(x)$ ($x\in X$) uniformly and independently from $\C{Y}{s}$.
Then
\[
\pr(\mbox{$L\cap F'$ admits an $X$-p.m.})\leq \pr(\mbox{$L\cap F$ admits an $X$-p.m.}) .
\]
\end{lemma}

\nin
(Here $s$ is arbitrary.  The lemma is also true, with essentially the same proof, 
if we replace $L\cap F$ and $L\cap F'$ by the usual ``$p$-random'' subgraphs of $F$, $F'$
gotten 
by keeping edges independently with probability $p$.)

Repeated application of Lemma~\ref{Lswitch} yields the statement we need:
\begin{cor}\label{Corollary}  For $X,Y,L$ as in Lemma~\ref{Lswitch}
and any $(d_x:x\in X)$, the probability that $L\cap F$ admits an $X$-p.m. is minimized 
over $F$ with $d_F(x)=d_x $ $\forall x\in X$
when, for some ordering of $Y$, each $x\in X$ is adjacent to the first $d_x$ vertices of $Y$.
\end{cor}

We return to $B'$ and the proof of \eqref{notext'}.
Let $|C|=\DD+x$ ($x\in \mathbb{Z}$) and observe:
\beq{dHbar}
d_{\ov{B'}}(v) ~\leq ~d^*_v~:=~|\nabla_G(v,V\sm C)|  ~=~ d_{H}(v)-x
\,\,\,\,\,\,\forall  v\in  C
\enq
(we just need $v\in C'$) and
\beq{sumvC'}
\sum_{v\in C'} d_v^* ~=~ \sum_{v\in C'} (d_{H}(v)-x)
~\le ~\sum_{v\in C} d_H(v)- |C'|x ~=~ 2\gz D^2 - |C'|x.
\enq

\mn

When $C$ is large (so $x>0$ and $\gz$ is large), 
\eqref{notext'}
is given by case (b) of Lemma~\ref{LMg} for $B=B'$,
with the hypotheses $J\approx D$, 
$d_{\oB}(v) =O(\eps D)$, $|\oB| \leq 2\gz D^2$ and (b) given
by (resp.)
\eqref{J'D}, [\eqref{dHbar} and \eqref{EDepsD}], [\eqref{dHbar} and \eqref{sumvC'} 
(the sacrifice of the extra $|C'|x$ in \eqref{sumvC'} being affordable here)], and 
the last line of \eqref{R'}.

Suppose instead that $C$ is small and let $\ell =-x$ ($\geq 0$).
Here it is enough to prove \eqref{notext'} assuming we have equality in 
\eqref{dHbar} (since \eqref{notext'} becomes ``less true'' if we delete edges from $B'$),
in which case we have 
$d_{\ov{B'}}(v)=d_H(v)-x\geq \ell$ ($\forall v\in C'$).

According to Corollary~\ref{Corollary}, it is then enough to show \eqref{notext'} under the 
additional assumption that there is some $\gG^*\in \C{\gG'}{\ell}$
all of whose vertices are isolated in $B'$.  
But in this case $B:=B'-\gG^*$ satisfies the hypotheses of Lemma~\ref{LMg}:
for $J\approx D$, 
$d_{\oB}(v) =O(\eps D)$ and $|\oB| \leq 2\gz D^2$, we again use 
\eqref{J'D}, [\eqref{dHbar} and \eqref{EDepsD}],
and [\eqref{dHbar} and \eqref{sumvC'},
in which---crucially---deleting $\gG^*$
removes the terms involving $x$];
and for (a) or (b) (as appropriate), we have 
$|X|=\DD-\ell = |\gG\sm \gG^*|=|Y|$ if $\gz$ is small, while for $\gz$ large 
\eqref{R'} gives $R \ge \eta D$.

This completes the proofs of \eqref{notext'} and \eqref{notext}.
\qed

\section{Proof of Lemma~\ref{Lsucceed}}\label{PLsucceed}

Recall the goal here is to show, for the ``process'' of Section~\ref{Clusters},
that each of (S1)-(S3) (defined just before Lemma~\ref{Lsucceed})
holds with probability $1-o(D/n)$.  The first of these (with (S1) repeated as \eqref{Lvgc1}) is easy:
\begin{obs}\label{ObsLvgc}
With probability $1-o(1/n)$,
\beq{Lvgc1}
|L_v\cap \{\gc_1\dots \gc_{m-1}\}| \leq 0.1 \gd \log n \,\,\, \forall v\in C.
\enq
\end{obs}
\nin
\emph{Proof.} 
This is an instance of Theorem~\ref{Cher'}, which,
with $\mu = m\kkk/\DD$ and $K\mu =0.1\gd\log n=:s$, gives (for any $v\in C$)
\[
\pr( |L_v\cap \{\gc_1\dots \gc_m\}| > s) 
~< ~\left(\frac{em\kkk}{Ds}\right)^s  ~=~ o(1/n),
\]
were we use
$m\ll D$ (see \eqref{mllD}) for the second inequality 
(so the actual bound is much better than $o(1/n)$).

\qed

Set
\beq{Qi}
Q_i=\{|L_v\cap \{\gc_1\dots \gc_{i-1}\}| < 0.1 \gd \log n \,\, \forall v\in C\};
\enq
thus $Q_1\supseteq \cdots \supseteq Q_m$ and $Q_m$ is the event in \eqref{Lvgc1}.
Setting 
$p = (\kkk -0.1\gd\log n)/D = (1+0.9\gd)\log n/D$,
we have
\beq{prvinJ}
\pr(v\in J_i|Q_i) > p \,\,\,\,\,\forall v\in C_{i-1}\sm N(\gs^{-1}(\gc_i)).
\enq
(The reason for the seemingly superfluous $Q_i$'s is that conditioning as in \eqref{prvinJ}
is problematic if we replace $Q_i$ by $Q_m $.)

\mn

The next lemma serves two purposes:  first, it says that (S2) holds 
with probability $1-o(D/n)$; second, it
will allow us to assume that, for $i$ with $|M_{i-1}|< \eta D$, 
$H_{i-1}$ is large,
which 
will be important when, in Proposition~\ref{PJanson}, we want to show that for each such $i$, $H[J_i]$
has a decent chance of being nonempty.

\begin{prop}\label{Pzi}
With probability $1-o(D/n)$, 

\nin
{\rm (a)}  $J_i\neq \0$ $\forall i\in [m]$ (that is, (S2) holds) and

\nin
{\rm (b)}
each
$z_i$ has $d_H(z_i) < 6\gz D/\gd$.
\end{prop}
\nin
\emph{Note.}  In (b) and the rest of this section, any assertion regarding a $z_i$ is 
is taken to include existence of $z_i$.

\mn
\emph{Proof.}
With
\[
E_i :=\{J_i=\0\}\vee\{d_H(z_i)\ge 6\gz D/\gd\}
\]
($i\in [m]$), we want
\beq{prei}
\pr(\cup E_i) =o(D/n).
\enq

From Observation~\ref{ObsLvgc} we have (recalling that $Q_m$ is the event in \eqref{Lvgc1}
and is contained in each $Q_i$)
\beq{prcupEi}
\pr(\cup E_i) \leq \pr(\ov{Q}_m)
+ \sum\pr (E_i\wedge Q_i) ~<~ o(1/n)  + \sum\pr(E_i|Q_i);
\enq
so, since $m\ll D$, \eqref{prei} will follow from (e.g.)
\beq{preiqi}
\pr(E_i|Q_i) < \exp[-(1+\gd/3)\log n].
\enq

\mn

Set
\[
Z=\{z\in C_{i-1}: z\not\sim \gs^{-1}(\gc_i), \,d_H(z) < 6\gz D/\gd\}.
\]
Since
\beq{zinC}
|\{z\in C: d_H(z) > 6\gz D/\gd\}| < \frac{2\gz D^2}{6\gz D/\gd} =\gd D/3,
\enq
we have
\[
|Z|~>~(1-b\eps)D - 2m- O(\eps D) -\gd D/3 ~>~ (1-\gd/2)D,
\]
where in the middle:
the first two terms lower bound $|C_{i-1}|$; 
the third is from Observation~\ref{Obsf}
(which implies $|\{v\in C: v\sim \gs^{-1}(\gc_i)\}| =O(\eps D)$); and the last is
given by \eqref{zinC}.  

Then for \eqref{preiqi} we just note that $E_i = \{J_i\cap Z=\0\}$
and apply \eqref{prvinJ}, yielding
\[
\pr(E_i|Q_i) ~=~ \pr(J_i\cap Z=\0|Q_i) 
~<~ (1-p)^{|Z|}
~<~ \exp[-(1+\gd/3)\log n] .
\]
\qed

Now define (for $i\in [m]$)
\[
R_i=\{\mbox{$d_H(z_j) < 6\gz D/\gd \,\,\, \forall z_j$ with $j<i$}\},
\]
noting that Proposition~\ref{Pzi}(b) gives 
\beq{printR}
\pr(\cap R_i) =1-o(D/n).
\enq
As was true for the $Q_i$'s (see \eqref{prcupEi}), \eqref{printR}
will allow us to focus on choices of $J_i$'s under conditioning on $R_i$'s
(as well as $Q_i$'s), so that the next assertion becomes the main point.

\begin{prop}\label{PJanson}
For each $i\in [m]$, if $|M_{i-1}|<\eta D$, then
\beq{prHJi}
\pr(H[J_i]=\0|Q_iR_i) < \exp\left[-\frac{\gz \kkk}{8b\eps}\right]  
\,\,\,\mbox{($= 1-q$; see \eqref{q}\emph{)}}.
\enq
\end{prop}

\nin
(As in (b) of ``Perspective'' following Lemma~\ref{Lsucceed},
this is roughly the best one can expect:  our assumptions allow
$H$ to be covered by some $T\sub C$ of size $O(\gz D^2/(\eps D)) =O(\gz D/\eps)$,
and the probability in \eqref{prHJi} is then at least about
$(1-k/D)^{|T|} =\exp[-O(\gz k/\eps)]$.)

\nin
\emph{Proof.}
Under conditioning on any $J_1\dots J_{i-1}$ for which  $Q_i R_i$ holds, we have
\[
\pr(v\in J_i) > p\,\,\,\,\,\forall
v\in C_{i-1}\sm N(\gs^{-1}(\gc_i))
\]
(see \eqref{prvinJ}) and 
\[
|H_{i-1}-N(\gs^{-1}(\gc_i))| \geq M:= \gz D^2 -[f(\gc_i)
+ \eta D\cdot 2b\eps D + m\cdot 6\gz D/\gd] \approx \gz D^2.
\]
Here the subtracted terms bound  
the number of edges $e$ in $H\sm (H_{i-1}-N(\gs^{-1}(\gc_i)))$: 
the first corresponds to $e$'s meeting $N(\gs^{-1}(\gc_i))$; the second, which uses 
$|M_{i-1}|<\eta D$ and \eqref{EDepsD}, to $e$'s meeting members of $M_{i-1}$; and 
the third, which uses $R_i$, to $e$'s meeting $z_j$'s with $j<i$.
The ``$\approx$'' is given by Observation~\ref{Obsf},
\eqref{etazeta} and \eqref{mllD}.

So it is enough to show 
\[
\pr^*(H[J_i]=\0) < \exp\left[-\frac{\gz \kkk}{8b\eps}\right],
\]
where we temporarily write $\pr^*$ for probabilities conditioned on (the worst case)
\[
\{|H_{i-1}-N(\gs^{-1}(\gc_i))| =M\}\wedge
\{\pr(v\in J_i)=p\,\,\,\forall v\in C_{i-1}\sm N(\gs^{-1}(\gc_i))\}.  
\]
This is an application of Theorem~\ref{TJanson},
for which we have (with some fairly arbitrary constants)
\[
\mu  = Mp^2> 2\gz \kkk^2/3
\]
and (since degrees in $H_{i-1}$ are less than $b\eps D$, again by \eqref{EDepsD}) 
\beq{ovgD}
\ov{\gD} < 2Mb\eps Dp^3 +\mu 
< 3b\eps \gz \kkk^3.
\enq
So the theorem gives
\[
\pr^*(H[J_i]=\0)  \le
\exp[-\mu^2/\ov{\gD}] < \exp\left[-\frac{\gz \kkk}{8b\eps}\right].
\]
\qed

It remains to show 
\beq{S3}
\pr(\mbox{(S3)}) = 1-o(D/n).
\enq
Let $\psi_1\dots \psi_m$ be independent ${\rm Ber}(q)$ r.v.s. and $\psi =\sum\psi_i$. 
We first show that \eqref{S3} will follow from
\beq{prPsi}
\pr(\psi < \eta D) < o(1/n).
\enq

To see this, set (for $i\in [m]$)
\[
\xi_i=\left\{\begin{array}{ll} 
0&\mbox{if $Q_iR_i$ fails or $|M_{i-1}|\geq \eta D$,}\\
\mbone_{\{H[J_i]=\0\}}&\mbox{otherwise,}
\end{array}\right.
\]
and notice that 
\[
\{\mbox{(S3) fails}\} ~= ~ \{|M_m|<\eta D\} 
~\sub ~ \cup \ov{Q_iR_i} ~\cup ~\left\{\sum \xi_i > m-\eta D\right\}
\]

\nin
(since $\{|M_m|<\eta D\} \cap\cap( Q_iR_i)$ implies $\sum\xi_i = m-|M_m|$).

On the other hand
Proposition~\ref{PJanson} says we can couple $(\xi_i)$ and $(\psi_i)$ with $\xi_i\leq 1-\psi_i$,
whence (using Observation~\ref{ObsLvgc}, as in \eqref{prcupEi}, together with
\eqref{printR})
\[
\pr(\mbox{(S3) fails}) ~\leq ~ \pr(\cup \ov{Q_iR_i}) + \pr (\sum\xi_i > m-\eta D)
~<~ o(1/n) + \pr (\psi<\eta D).
\]\qed

\nin
\emph{Proof of \eqref{prPsi}.}
We have $\E\psi =mq =K\eta D$ (see \eqref{KetaDq}), so Theorem~\ref{T2.1} gives (say)
\[
\pr(\psi <\eta D) < \exp[-K\eta D/3].
\]
This gives \eqref{prPsi} if (say) $\gz k/\eps \leq 1$, since in this case $q =\Theta(\gz k/\eps)$ and,
by \eqref{K},
\[
K\eta D\gg \eps q/\gz = \Theta (k) \,\,\,\,(=\Theta(\log n)).
\]

For $\gz k/\eps > 1$, we need a better (naive) bound:  
setting $\vt = 1-q =\exp\left[-\frac{\gz \kkk}{8b\eps}\right]$,
and using $m=K\eta D/q$,
we have
\[
\pr(\psi < \eta D) < \C{m}{\eta D} \vt^{m-\eta D} 
<
\exp\left[\eta D\log (eK/q) -\Theta \left((\gz k/\eps)K\eta D\right)\right]
=o(1/n).
\]
Here the last inequality (an understatement) is again based on \eqref{K}, which in ``exp'' gives:

\nin
(a)  $K\eta D\gg \eps q/\gz =\gO(1)$ (since $\gz=O(\eps)$ in general, and in the 
present regime $q=\gO(1)$), so the second term, say $T$, satisfies $T\gg k$; and

\nin
(b)  $T$ dominates the first term, since $K\gg 1$ (by \eqref{K})
and $\gz k/\eps>1$, which also implies $q=\gO(1)$.

\qed

This completes the proof of Lemma~\ref{Lsucceed}.

\section{Matchings}\label{Matchings}

\mn
\emph{Proof of Lemma~\ref{LMg}.}
As usual for such statements,
we want to show $K$ is likely to satisfy Hall's Condition, which we reformulate as:
\beq{HCA}
\QQQ\sub X, |\QQQ|\leq \lceil J/2\rceil ~\Ra ~ |N_K(\QQQ)|\geq |\QQQ|
\enq
and
\beq{HCB}
\QQQ\sub Y, ~1+R\leq |\QQQ|\leq \lfloor J/2\rfloor +R~\Ra ~ |N_K(\QQQ)|\geq |\QQQ|-R.
\enq

The event in \eqref{HCA} fails at $\QQQ\in\C{X}{\uuu}$ 
iff there is $Z\in\C{Y}{\uuu-1}$ for which
\beq{LvY}
M_v\sub Z\cup N_{\ov{B}}(v) \,\,\,\,\,\,\forall v\in \QQQ,  
\enq
the probability of which is less than
\[
\C{J+R}{\uuu-1}\prod_{v\in \QQQ} \left(\frac{\uuu-1 + d_{\ov{B}}(v)}{J+R}\right)^t
<
\C{J+R}{\uuu-1}\left(\frac{\uuu + O(\eps D)}{J+R}\right)^{t\uuu}.
\]
So the overall probability that \eqref{HCA} fails is less than 
\beq{SumJ2}
\sum_{\uuu=1}^{\lceil J/2\rceil}\C{J}{\uuu}\C{J+R}{\uuu-1}
\left(\frac{\uuu + O(\eps D)}{J+R}\right)^{t\uuu},
\enq
which is easily seen to be (conservatively) $o(1/n)$, the summand being largest when $\uuu=1$.

For $\QQQ$'s contained in $Y$ we observe that if $K$ (which for this observation could be 
any bigraph on $X\cup Y$) violates \eqref{HCB}, and 
$I=\{y\in Y: \mbox{$y$ is isolated in $K$}\}$,
then for any $a\in [0,R+1]$ either
\beq{either}
|I|\geq a\enq
or, for some $u\in [2,J/2]$,
\beq{or}
\mbox{$\exists$ $A\in \C{Y\sm I}{u+R-a+1}~$ 
with $~|N_K(A)|\leq u-1$.}
\enq
[If $A'\in \C{Y}{u+R}$ violates \eqref{HCB} (so $u\leq J/2$) and has fewer than $a$ 
isolates (so failure of \eqref{HCB} implies 
$u\geq 2$), then any $A\in \C{ A'\sm I}{u+R-a+1}$ is as in \eqref{or}.]

Here, if we are in (a) (of Lemma~\ref{LMg}, so with $R=0$), we take $a=1$
(which puts us in essentially the standard argument for a p.m.\ in a random bigraph),
and if in (b) take any $a$ with
\beq{takeany}
\max\{\gz D,1\}\ll a\ll \eta D \,\,\, (\lesser R)
\enq
(this is possible by \eqref{etazeta}),
and should show that each of \eqref{either}, \eqref{or} is unlikely.
We now (re)define 
\[
q = t/|Y| = (1+\gd')\log n/(J+R)
\]
($= \pr(y\in M_x)$ for any $(x,y)\in X\times Y$).

\mn
\emph{For} \eqref{either}:  
The probability that all $y\in Z\in \C{Y}{a}$ are isolated is less than 
$
(1-q)^{|\nabla_B(X,Z)|} \leq (1-q)^{Ja-|\ov{B}|},
$
so
\[   
\pr(\eqref{either}) ~\leq ~ \C{J+R}{a}(1-q)^{Ja-|\ov{B}|} < \left[2J n^{-(1+\gd'/2)}\right]^a
=o(D/n).
\]   
Here the second inequality uses $|\ov{B}|\ll Ja$, which holds in (a) (where $a=1$) since 
$\gz D\ll 1$, and in (b) since $|\ov{B}| \leq 2 \gz D^2$ and, by \eqref{takeany},
$a\gg \gz D$.

\mn
\emph{For} \eqref{or}:
We have
\beq{Por}
\pr(\eqref{or})  < 
\Cc{J+R}{u+R-a+1}\Cc{J}{u-1}(1-q)^{(u+R-a+1)(J-u+1)-2\gz D^2} (uq)^{u+R-a+1}.
\enq
[The first factor corresponds to (overestimates) possibilities for $A $ 
and the second to possibilities for $Z\in \C{X}{u-1}$ containing $N_K(A)$;
the third bounds $\pr(\nabla_K(X\sm Z, A)=\0)$, and the fourth the probability that 
each $y\in A$ has a neighbor in $Z$.  (Note these last two events are independent.)]

We may bound the second factor in \eqref{Por} by $(eJ/(u-1))^{u-1}$
and the third by 
\beq{thirdbd}
n^{-(1+\gd'/2)(u+R)(1-u/J)}.
\enq
[For \eqref{thirdbd} we have
\[   
(u+R-a+1)(J-u+1)-2\gz D^2 \approx (u+R)(J-u),
\]   
as follows from: (i) $u+R-a+1\approx u+R$
(for (a) this is just $u=u$, and for (b) it holds because
$a\ll R$; see \eqref{takeany}); (ii) $J-u=\Theta(D)$ (since $u\leq D/2$); and 
(iii) $\gz D^2\ll u+RD$ (by \eqref{epszeta} in (a) 
and \eqref{takeany} in (b)).
A second application of $R\ll D$
then gives \eqref{thirdbd}.]

So the product of these two factors is at most (say)
\beq{2bds}
\begin{array}{ll} 
D^{(u-1)} n^{-(1+\gd'/5)(u+ R)} &\mbox{if $u< \gd' J/4$,}\\
n^{-u/2}&\mbox{otherwise.}
\end{array}
\enq
On the other hand, the product of the first and last factors is less than
\[
\left[\tfrac{e(J+R)}{u+R-a+1}\cdot \tfrac{u(1+\gd') \log n}{J+R}\right]^{u+R-a+1}
= O(\log n)^{u+R},
\]
which is negligible relative to the bounds in \eqref{2bds}.
In particular, $\pr(\eqref{or}) =o(1/n)$.

\qed

\mn
\emph{Proof of Lemma~\ref{Lswitch}.}
All graphs here are bipartite on $X\cup Y$.  
Lemma~\ref{Lswitch} is proved by coupling---which is as one should expect, 
but does require some care.  

Say a bigraph $J$ (on $X\cup Y$) is \emph{good} if it admits an $X$-p.m.
For any $J$ and $Z\sub X$, let $J(Z)$ be gotten from $J$ by ``switching'' $Z$, meaning
$J$ and $J(Z) $ agree \emph{off} $\nabla(Z,\{u,v\})$ (i.e., on edges not meeting 
both $Z$ and $\{u,v\}$),
\[
\mbox{$N_{J(Z)}(v)\cap Z= N_{J}(u)\cap Z~$,}
\]
and similarly with $u$ and $v$ interchanged.
Notice that, for any $J$ and $Z$, 
\beq{couple1}
\pr(L=J) =\pr (L=J(Z)).
\enq

Let 
\[
\mbox{$D=N_F(u)\cap N_F(v)$, $A=N_F(u)\sm N_F(v)$, and $B=N_F(v)\sm N_F(u)$.}
\]
From this point $J$ will be a possible $L$ (i.e.\ value of the random $L$ of the lemma).  
We will couple the pairs $\{J,J(D)\}$ and $\{J(B),J(B\cup D)\}$.
Notice that each of these two collections of ``pairs'' (some of which---those with 
$N_J(u) \cap D= N_J(v) \cap D$---are 
actually singletons) partitions the set of possible $L$'s.
In view of \eqref{couple1},
it is thus enough to prove the next (deterministic) assertion, which is the main point here.
\begin{claim}\label{Ccouple}
For every $J$,

\nin
(i) if either of $J(B)\cap F'$, $J(B\cup D)\cap F'$ is good, then at least one of 
$J\cap F$, $J(D)\cap F$ is good;

\nin
(ii) if $J(B)\cap F'$ and $J(B\cup D)\cap F'$ are good, then 
$J\cap F$ and $J(D)\cap F$ are good.
\end{claim}

\nin
[For perspective:  switching $B$ is natural in this context, since 
edges of $\nabla_J(B,v)$, which can be useful when we intersect with $F$, 
become edges of $\nabla_{J(B)}(B, u)$, which can be useful in $J(B)\cap F'$. 
(One might guess it would 
be enough here to just couple $J$ and $J(B)$, but this doesn't work.)]

\nin
\emph{Proof of Claim~\ref{Ccouple}.}
For the rest of this discussion we fix $J$.
We first dispose of two trivial cases.  First, since $H:=Q-\{u,v\}$ is the same for 
all four graphs $Q$ in Claim~\ref{Ccouple}, we may assume $H$ is not good.
Second, if one of $J(B)\cap F'$, $J(B\cup D)\cap F'$ admits an $X$-p.m.\ using just one of 
$u,v$, say $M\cup \{xw\}$ with $M\sub H$, $x\in A\cup B\cup D$ and $w\in \{u,v\}$, then
$M$ is easily seen to extend to an $X$-p.m.\ using one of $xu$, $xv$ in each of $J\cap F$, $J(D)\cap F$.

So we may assume that any $X$-p.m.\ of either of $J(B)\cap F'$, $J(B\cup D)\cap F'$ uses both 
$u$ and $v$.
A little reformulation will now be helpful.
Say (with $H$ as above) $\{s,t\} \sub X$ (with $s\neq t$) is \emph{exposable} if 
\beq{exp} 
\mbox{$\KKK-\{s,t\}$ contains an ($X\sm\{s,t\})$-perfect matching;}
\enq
$(s,t)$ is exposable if $\{s,t\}$ is; and $(P,R)$ (with $P,R\sub X$) is 
exposable if $P\times R$ contains an exposable pair.
(It will be helpful to have these notions for ordered pairs, 
though the order doesn't matter.)

It is immediate that if $Q$ is any of the graphs in Claim~\ref{Ccouple}, then
\beq{Qexp}
\mbox{$Q$ is good iff $(N_Q(u), N_Q(v))$ is exposable.}
\enq

Let $N_J(u) =A_u\cup B_u\cup D_u$ (with $A_u\sub A$ and so on), and similarly
for $N_J(v)$, and, to slightly simplify notation, set 
\[
\mbox{$A_u=A', \, B_v=B',\, D_u=S,\,$ and $\, D_v=T$}
\]
($A_v$ and $B_u$ play no role here);
rename the graphs in Claim~\ref{Ccouple}
\beq{rename}
\mbox{$K_1=J\cap F, \, K_2= J(D)\cap F,\, K_1' =J(B)\cap F',\,$ and $\, K_2'= J(B\cup D)\cap F'$;}
\enq
and write $N_1$ for neighborhood in $K_1$ and so on.

Each of the graphs in \eqref{rename} 
agrees with $\KKK$ off $\{u,v\}$,
while for the neighborhoods of $u$ and $v$ we have
\[
N_1(u) = A'\cup S,\,  N_1(v) = B'\cup T, \, 
N_2(u) = A'\cup T,\,  N_2(v) = B'\cup S
\]
and
\[
N'_1(u) = A'\cup B'\cup S,\,  N'_1(v) = T, \, 
N'_2(u) = A'\cup B'\cup T,\,  N'_2(v) = S.
\]
So in view of \eqref{Qexp}, Claim~\ref{Ccouple} is equivalent to the next statement,
which is the one we will prove.
\begin{claim}\label{Ccouple'}
\nin
(i) If either of
$(A'\cup B'\cup S,\,  T),\, (A'\cup B'\cup T,\,  S)$ is exposable, then at least one of 
$(A'\cup S,B'\cup T),\, (A'\cup T,B'\cup S)$ is exposable.

\nin
(ii) If $(A'\cup B'\cup S,\,  T)$ and $(A'\cup B'\cup T,\,  S)$ are exposable, then 
$(A'\cup S,B'\cup T)$ and $( A'\cup T,B'\cup S)$ are exposable.
\end{claim}

\nin
\emph{Proof of Claim~\ref{Ccouple'}.}
First observe that if $(x,y)\in (A'\cup B'\cup S)\times T$ is exposable, then 

\nin
(a) $(A'\cup S,B'\cup  T)$ is exposable \emph{unless} $x\in B'$ and $y\not\in S$

\nin
(e.g.\ if $x\in B'$ and $y\in S$, then we use $(y,x)\in (A'\cup S)\times (B'\cup T)$), and

\nin
(b) $(A'\cup T,B'\cup  S)$ is exposable \emph{unless} $x\in A'$ and $y\not\in S$.

This and the corresponding statement for an exposable $(w,z)\in (A'\cup B'\cup T)\times S$
immediately give (i), and they also give (ii)---where we are given \emph{both} $(x,y)$ and $(w,z)$
as above---\emph{unless} either 
\beq{xyzw1}
\mbox{$x\in A',\, y\in T\sm S,\, w\in B'\,$ and $\, z\in S\sm T$}
\enq
(where the observations at (a) and (b) only say that each of $\{x,y\}$, $\{w,z\}$
gives exposability of $(A'\cup S,B'\cup T)$) or 
\[   
\mbox{$x\in B',\, y\in T\sm S,\, w\in A'\,$ and $\, z\in S\sm T$.}
\]   
Assume (w.l.o.g.) that \eqref{xyzw1} holds---so $(A'\cup S,B'\cup  T)$ is exposable---and
observe (check) that exposability of any one of the four remaining pairs (those
other than $\{x,y\}, \{w,z\}$) from $\{x,y,w,z\}$ will give exposability of $(A'\cup T,B'\cup  S)$.

Here, finally, we do need to slightly think about structure.
We are given matchings $M_1,M_2$ witnessing exposability of 
$\{x,y\}$ and $\{w,z\}$ (e.g.\ $M_1$ is an ($X\sm\{x,y\}$)-p.m.\ of $\KKK-\{x,y\}$) ,and
consider the alternating paths (in $M_1\cup M_2$) using $x,y,w,z$.
If one of these paths, $P$, joins two of $x,y,w,z$---necessarily one from each 
of $\{x,y\}, \{w,z\}$, so say $x$ and $w$---then (e.g.) $M_1\Delta P$ witnesses exposability of 
$\{w,y\}$.  And if not then we get the same result using $M_1\Delta (P_x\cup P_w)$
(where $P_x$, $P_w$ are the paths containing $x$ and $w$, each of which has its other end in $Y$).
So in either case 
we have the desired additional exposable pair.

This completes the proofs of Claim~\ref{Ccouple'}, Claim~\ref{Ccouple}, 
\and Lemma~\ref{Lswitch}.
\qed

\section{Discussion}  \label{Disc}

Apart from its algorithmic aspect, ``palette sparsification'' (though not under this name)
has been of independent interest at least since 
work of Krivelevich and Nachmias in \cite{KN1,KN2}; see 
\cite{CH} for a few references.
We briefly mention two strands, the first just to point out the connection, the second 
closer to what we do here.

The following important result of Keevash \cite{Keevash}
and Jain and Pham \cite{JP}, improving 
work of Sah, Sawhney and Simkin \cite{SSS} and
Kang, Kelly, Kuhn, Methuku and Osthus \cite{DKKMO}, had been conjectured by 
Casselgren and H\"{a}ggkvist \cite{CH} and Luria and Simkin \cite{LS}, and is
the natural answer to a question of Johansson \cite{Joh}.

\begin{thm}\label{TKJP}
There is a fixed C such that if $L_e$'s, $e\in E(K_{n,n})$, are drawn uniformly and independently
from the ($C\log n$)-subsets of
$\gG=[n]$, then G is L-\textbf{edge}-colorable w.h.p.
\end{thm}
\nin
The papers above prove more general results, regarding thresholds for 
Steiner triple systems, Latin squares (addressed by Theorem~\ref{TKJP})
and more.  See e.g.\ the discussion in \cite[Sec.\ 1.1]{DKKMO} for many intriguing, still-open 
possibilities, including asymptotic (even hitting time) statements akin to Theorem~\ref{TKK}.  
(The approach of \cite{SSS,DKKMO,Keevash,JP}, based on \cite{FKNP,PP}, is 
unrelated to anything in the present work or its antecedents
\cite{ACK,AA}.)

Closer to present concerns, 
Alon and Assadi \cite{AA} consider the general (and from a list-coloring standpoint very natural)
situation in which $\gG$ is arbitrary and each $L_v$ is a uniform $\ell$-subset of some $S_v\sub \gG$.
They show (with $G$ again on $V=[n]$ and $d_v:=d_G(V)$):
\begin{thm}\label{TAA}
In the above framework, $G$ is L-colorable w.h.p.\ if either

\nin
{\rm (a)}  $|S_v|=(1+\eps)d_v$ ($\forall v$) and $\ell = C_\eps\log n$ 
(for any fixed $\eps >0$ and a suitable $C_\eps$), or

\nin
{\rm (b)}  $S_v=[d_v+1]$ ($\forall v$) and $\ell = C\log n$ 
(for a suitable $C$).
\end{thm}

At this writing we \emph{believe} we can extend Theorem~\ref{TKK} to the above framework:

\begin{conj}\label{CKK}
For fixed $\eps>0$, 
if 
$\gD_G \leq D$, $|S_v|=D+1$ $\forall v$, and $\ell =(1+\eps)\log n$, 
then $G$ is L-colorable w.h.p.
\end{conj}

\nin
We also suspect that this extends even to the ``correct'' version of Theorem~\ref{TAA}
(eliminating the $\eps$ in (a) and the somewhat arbitrary choice of $S_v$'s in (b)):

\begin{conj}\label{CKK2}
For fixed $\eps >0$,
if $|S_v|=d_v+1$ $\forall v\in V$ and $\ell = (1+\eps)\log n$,
then G is L-colorable w.h.p.
\end{conj}

\end{document}